\documentclass{amsart}
\usepackage{amssymb}

\usepackage[%
 colorlinks=true,%
 linkcolor=blue,%
 citecolor=blue,%
 filecolor=blue,%
 urlcolor=blue,%
]{hyperref}

\usepackage{graphics}

\usepackage[
nonotes
]{degt}

\let\=\B

\def\MAT{\operatorname{Mat}_{2\times2}}

\def\paragraph{\subsubsection{}}
\def\loccit#1{(\ref{#1})}
\def\iref#1#2{\ref{#1}\loccit{#2}}
\let\ditto\loccit

\let\minusplus\mp

\def\Prod{\BG3\cdot\Inn\AA}
\def\SS{\Bbb S}

\def\rdeg{\operatorname{wd}}  
\def\bdeg{\operatorname{dg}}  

\let\Ga\alpha
\let\Gb\beta
\let\Gg\gamma
\let\Gd\delta
\let\Gl\lambda
\let\Gf\varphi
\let\Gs\sigma
\let\Gr\rho

\def\scirc{^\circ}

\def\tj{\tilde\jmath}

\def\be{\bold e}
\def\bv{\bold v}

\def\CM{\Cal M}

\def\CR{\Cal R}

\def\X{\Bbb X}
\def\Y{\Bbb Y}

\def\CV{\Cal V}
\def\CU{\Cal U}
\def\IM{\bar\CV}
\def\tIM{\CV}

\let\BB=B  
\let\CC=C  
\let\EE=E  
\let\DD=D  
\let\disk\Delta
\let\slope\varkappa

\def\BM{\frak{Im}} 
\let\AM=A          
\def\AM{\mathrm{A}} 
\def\bAM{\bar\AM}
\def\bDelta{\bar\Delta}

\def\GDG{\Bbb D}         
\def\FF#1{\Bbb F_{#1}}   
\def\AA{\Bbb G}          
\def\AA{\frak F}         
\let\cp\Phi              

\def\bm{\frak m}

\def\sproj{\suptext{proj}}
\def\saff{\suptext{afn}}
\def\piproj#1{\pi\sproj_{#1}}
\def\piaff#1{\pi\saff_{#1}}

\def\inserthyphen{\ifcat\next a-\fi\ignorespaces}
\let\BLACK\bullet
\let\WHITE\circ
\def\CROSS{\vcenter{\hbox{$\scriptstyle\mathord\times$}}}
\let\STAR*
\def\TRIANG{\vcenter{\hbox{$\scriptstyle\mathord\vartriangle$}}}
\def\pblack-{$\BLACK$\futurelet\next\inserthyphen}
\def\pwhite-{$\WHITE$\futurelet\next\inserthyphen}
\def\pcross-{$\CROSS$\futurelet\next\inserthyphen}
\def\pstar-{$\STAR$\futurelet\next\inserthyphen}
\def\ptriang-{$\TRIANG$\futurelet\next\inserthyphen}
\def\black{\protect\pblack}
\def\white{\protect\pwhite}

\def\NO#1{n_{#1}}
\def\nblack{\NO\BLACK}
\def\nwhite{\NO\WHITE}

\def\Mat#1{\bmatrix#1\endbmatrix}
\def\mat[#1,#2,#3,#4]{\Mat{#1&#2\\#3&#4}}
\def\tf{\tilde\varphi}
\let\Bs\beta
\def\Bs{\bar\Gs}
\def\Bb{\bar\Gb}
\def\bG{\bar G}
\def\plus{^+}
\def\minus{^-}
\def\two{^{\text{bu}}}
\let\mp\varkappa  

\let\vI\I
\let\vII\II
\let\vIII\III
\let\vIV\IV
\def\vIIex{\vII_{\fam0 ex}}
\let\vZ0
\def\kI{k^{\vI}}
\def\kIII{k^{\vIII}}
\def\kIV{k^{\vIV}}
\def\kall{k^{\operatorname{all}}}
\def\KI{K^{\vI}}
\def\KIII{K^{\vIII}}
\def\KIV{K^{\vIV}}

\def\ee{\operatorname{e}}

\def\reg(#1){(\!(#1)\!)}
\def\dR{\partial R}

\def\bbbk{\smash{\tilde\Bbbk}}
\def\bbbk{\Bbb K}

\title{The Alexander module of a trigonal curve}

\author{Alex Degtyarev}

\address{%
Bilkent University\\
Department of Mathematics\\
06800 Ankara, Turkey%
}

\email{degt@fen.bilkent.edu.tr}

\keywords{%
Trigonal curve,
fundamental group,
Alexander module,
Alexander polynomial,
Burau representation,
modular group%
}

\subjclass{%
Primary: 14H30; 
Secondary: 14H45, 
14H50, 
20F36
}

\begin{document}

\begin{abstract}
We describe the Alexander modules and Alexander polynomials (both
over~$\Q$ and over finite fields~$\FF{p}$) of generalized trigonal
curves. The rational case is closed completely; in the case of
characteristic $p>0$, a few points remain open.
The results obtained
apply as well to plane curves with deep singularities.
\end{abstract}

\maketitle

\section{Introduction\label{S.intro}}

\subsection{Motivation}\mnote{renamed, extended,
motivation added; abstract extended}
This paper continues the
systematic study of the fundamental groups of
(generalized) trigonal curves
that was started in~\cite{degt:trigonal}.
(By the common abuse of the language, speaking about the fundamental group of
an embedded curve,
one refers to the group of
the complement of the curve;
see Subsection~\ref{s.pi1} for the precise description of the groups
to be studied.)
The principal motivation for this research is the belief that
there should be
strong restrictions to the complexity of
these groups, far beyond the obvious fact
that they admit presentations with at most three
generators.
Thus, only about a dozen of distinct groups appear as the fundamental groups
of irreducible plane sextics with a triple point
(see~\cite{degt:dessin} and references therein), which are a special class of
generalized trigonal curves.
(Remarkably, the commutants of most finite groups
obtained in this way are of the form $\SL(2,\Bbbk)$, where $\Bbbk$ is a
finite field.)
These restrictions are due to the fact that the monodromy group of a trigonal
curve is a genus zero subgroup of the modular group, see
Subsection~\ref{s.sections} and \autoref{th.monodromy}; hence, it is
sufficiently `large', resulting in a sufficiently small fundamental group.
At present, it is not quite clear how or even in what terms such
fundamental groups can
be characterized; as a first step, we make an attempt to describe their
metabelian invariants.

Another\mnote{rewritten}
special feature of trigonal curves is the fact that,
in this case,
the relation between
the fundamental group and the geometry of a curve
is `two-sided', as all curves
with `at least' a certain fundamental group are essentially
induced from some universal curve with this property,
see Speculation~1.2.1 and a number of examples
in~\cite{degt:trigonal}.
For example,\mnote{example added}
\cite[Theorem 1.2.5]{degt:trigonal} characterizes the so-called
curves of \emph{torus type} in terms of their Alexander polynomial;
remarkably, a very similar
assertion holds for irreducible plane sextics, see~\cite{degt:Oka}.
An essential intermediate statement concerning the universal curves
is cited in
\autoref{th.monodromy}.

A\mnote{rewritten}
generalized trigonal curve in the Hirzebruch surface~$\Sigma_1$
(plane blown up at one point)
can be
regarded as a curve in the plane $\Cp2=\Sigma_1/E$, where $E$ is the
exceptional section, and as such it has a distinguished singular point of
multiplicity $(\text{degree}-3)$, see Subsection~\ref{s.plane}.
Thus, the study of trigonal curves sheds light to
the classical problem about the fundamental group of a plane
curve. (It is this construction that motivated my original interest in
trigonal curves.)
As an example,
the passage to the trigonal
model, combined with the techniques of \emph{dessins d'enfants}
described below, lets one
compute the fundamental groups of all irreducible sextics with a singular
point of multiplicity at least three, see~\cite{degt:dessin}, whereas the
groups of a number of sextics with double singular points only are
still unknown.
It is worth mentioning that
there is a mysterious similarity, although not quite literal
coincidence, between the properties of plane sextics and those of trigonal
curves (see~\cite{degt:trigonal} for a more detailed discussion);
it must be due to the similarity between $K3$- and elliptic surfaces.

The principal tool\mnote{edited} used in the paper is the
correspondence between trigonal curves in Hirzebruch
surfaces, genus zero subgroups of the modular group, and a certain
class of planar bipartite ribbon graphs (essentially, Grothendieck's
\emph{dessins d'enfants} of the modular $j$-invariant), see,
\eg, \cite{Birch,Bogomolov,degt:trigonal,degt:monodromy,Kulkarni}.\mnote{refs
added}
As a by-product, we obtain some information on
the scarcity of the image of the
Burau representation of the braid group~$\BG3$, see
\autoref{rem.Burau} on the `Burau congruence subgroups',
although no attempt to formalize these
results has been made.

\subsection{The subject}\mnote{separate subsection}
In~\cite{degt:trigonal}, we gave a complete classification of
the dihedral quotients of
the fundamental group of a generalized trigonal curve.
Here, we deal with the ultimate metabelian
invariants of a curve, \viz. its so-called Alexander module and Alexander
polynomial. In the context of algebraic curves, this concept
appeared essentially in~\cite{Zariski:group}; it was later
developed
in~\cite{Libgober:Duke,Libgober:invariants,Libgober:modules,Libgober:Burau},
and it has been a subject of intensive
research since then, see recent surveys~\cite{Libgober:survey2,Oka:survey}
for further references.

For an irreducible
generalized trigonal curve $\CC$ in the Hirzebruch surface~$\Sigma_d$
(see Section~\ref{S.curves}),\mnote{name and ref added}
the
\emph{Alexander module} $\AM_\CC$ can be defined
as the homology group $H_1(X)$
of the maximal cyclic covering $X\to\Sigma_d$ ramified at~$\CC$
and the exceptional section~$E$, see Subsections~\ref{s.kernels}
and~\ref{s.pi1} for
details.
The deck translation
automorphism of the covering induces an action on~$\AM_\CC$,
turning it into a module over the ring $\Lambda:=\Z[t,t\1]$
of Laurent polynomials.
This module describes the fundamental group
of the curve modulo its second commutant.
Classically, one tensors~$\AM_\CC$ by~$\Q$ to
get
a torsion
module over the principal ideal domain $\Lambda\otimes\Q$; the
order~$\Delta_\CC$ of
$\AM_\CC\otimes\Q$
is called the \emph{Alexander
polynomial} of~$\CC$. To capture the
integral torsion
of~$\AM_\CC$, we
will also consider the product $\AM_\CC\otimes\FF{p}$ for a
prime~$p$; the order
$\Delta_{\CC,p}\in\Lambda\otimes\FF{p}$
of this product is
called the \emph{$({\bmod}\,p)$-Alexander polynomial}.
(A similar approach was used in~\cite{Libgober:invariants}, where some
$({\bmod}\,p)$-Alexander polynomials were computed.)

As\mnote{a very brief survey; refs added}
in the knot theory, the Alexander polynomial
is a purely algebraic invariant of the fundamental group of the curve,
but it is usually much easier to compute directly.
The
classical rational
polynomial~$\Delta_C(t)$ can be computed
by means of the Hodge theory,
in terms of the superabundance of certain linear
systems related to the singularities of the curve, see
\cite{degt:poly,Esnault,Libgober:invariants,Loeser.Vaquie}.
(Although most results
are stated for plane curves, they can easily be adapted to curves in any
surface.)
Besides, there are a great deal of the so-called \emph{divisibility
theorems}, bounding the Alexander polynomial in terms of the degree of the
curve and/or its singularities.
Some of these theorems, \eg,~\cite{Libgober:Duke,Libgober:invariants}, are of
purely topological nature and apply as well to pseudo-holomorphic curves and
$({\bmod}\,p)$-Alexander polynomials.
Others, \eg,~\cite{degt:divide}, rely upon the vanishing theorems in
algebraic geometry; they give better estimates, but work for algebraic curves
and rational Alexander polynomial only.
All these statements are in sharp contrast with the principal results of this
paper, as we show that, for each~$p$, the $({\bmod}\,p)$-Alexander polynomial
of a \emph{trigonal} curve may take but finitely many values, no matter what
the singularities are. The particular case $p=0$, see \autoref{th.0}, can
be translated into a certain restriction to the complexity of the
singularities of a trigonal curve and their mutual position: the
superabundance of some linear systems cannot be too large.

\subsection{Principal results}
Throughout the paper, we assume that $p$ is a prime or zero and let
$\Bbbk_0=\Q$ and $\Bbbk_p=\FF{p}$ for $p>0$. (When $p$ is fixed,
we abbreviate~$\Bbbk_p$ to~$\Bbbk$.)
For an
element~$\xi$ algebraic over~$\Bbbk$, we denote by
$\mp_\xi\in\Bbbk[t]$ its minimal polynomial and,
if $\xi$ is understood, we let
$\bbbk=\Bbbk(\xi)=(\Lambda\otimes\Bbbk)/\mp_\xi$.\mnote{$\bar\Bbbk\to\bbbk$
everywhere}
The cyclotomic polynomial (over~$\Q$) of order~$n$ is denoted
by~$\cp_n$.

As this paper is just a first step towards the understanding of
the Alexander module, we choose to work over a field and
consider the\mnote{italic to emphasize
that {\bf this} is the definition} \emph{specializations}
$\AM_\CC(\xi):=(\AM_\CC\otimes\Bbbk)/\mp_\xi$, see
Subsection~\ref{s.specialization},
thus reducing to $r=1$ higher torsion summands
of the form
$\CG{p^r}$ or
$(\Lambda\otimes\Bbbk)/\mp_\xi^r$, $r>1$, which may and do appear
when $p>0$.
In other words,\mnote{explanation added}
we are trying to enumerate the possible roots~$\xi$ of the Alexander
polynomial~$\Delta_{C,p}$ or, equivalently, its irreducible
factors, which are of the form~$\mp_\xi$.
Note\mnote{note on the dimension}
that $\AM_C(\xi)$ is a vector space over~$\bbbk$,
and therefore we can speak about
its dimension rather than rank.

\convention
Since\mnote{convention/explanation added; $N$ defined}
$\Delta_{C,p}$ is defined over~$\Bbbk_p$ itself, the set of its roots is
Galois invariant. For this reason, in most statements we refer to the minimal
polynomials $\mp_\xi\in\Bbbk_p[t]$ rather than to particular roots
$\xi\in\bbbk_p$. With~$\xi$ or $\mp_\xi$ understood, we fix the
notation~$N$ for the multiplicative order $\ord(-\xi)$.
Certainly, $N$ is determined by~$p$ and $\mp_\xi$; however, in view of the
importance of this parameter, we will speak about triples $(p,N,\mp_\xi)$
rather than just pairs $(p,\mp_\xi)$ (or even singletons~$\mp_\xi$, which
formally remember~$\Bbbk_p$ as their coefficient field).
It is worth mentioning that each pair $(p,N)$, $N\ge1$, corresponds to but
finitely many minimal polynomials $\mp_\xi$, \viz. the irreducible divisors
(over~$\Bbbk_p$) of $\cp_N(-t)$, and in some statements it is $(p,N)$ that is
fixed/discussed, whereas $\mp_\xi$ is allowed to vary.
\endconvention

The principal results of the paper are summarized in the next
four statements.
We close completely the case $p=0$, while for $p>0$ a
certain range still remains open.
Conjecturally, the Alexander polynomial of a non-isotrivial
trigonal curve can take finitely many values, and all irreducible
factors are indeed listed in the paper (with
\autoref{tab.examples} in \autoref{example.7--10}
taken into account).
Note that, unlike a number of known divisibility theorems
(\cf.~\cite{Libgober:Duke,Libgober:invariants,degt:divide}), the
bounds below are \emph{universal}, as we do not make any assumptions
about the singularities of the curve or its degree.

\theorem\label{th.0}
The Alexander polynomial~$\Delta_\CC$
of an irreducible non-isotrivial
generalized trigonal curve~$\CC$ can take only the following four
values\rom:
$\cp_6$, $\cp_6^2$, $\cp_{10}$, and~$\cp_{10}^2$. All four values
can be realized by genuine trigonal curves.
\endtheorem

\theorem\label{th.p}
Let $p>0$, and assume that the $({\bmod}\,p)$-Alexander
polynomial
$\Delta_{\CC,p}$ of a non-isotrivial
generalized trigonal curve~$\CC$ has
a root $\xi\in\bbbk\supset\Bbbk_p$.\mnote{definition of $N$ moved above}
Then, with the
exception of the
fourteen
triples $(p,N,\mp_\xi)$ listed in
\autoref{tab.factors}, one has $1\le N\le10$.
If $\CC$ is irreducible and
$N\ne3$ or~$5$,
one has $\dim_{\bbbk}\AM_\CC(\xi)=1$.\mnote{$\rank\to\dim$}
\table
\caption{Exceptional factors of~$\Delta$ ($N>10$)}\label{tab.factors}
\bgroup
\def\*{\llap{$^{*}$}}%
\let\comma,\catcode`\,\active\def,{$\comma\ \ $}%
\let\scolon;\catcode`\;\active\let;\comma
\def\get[#1]{$#1$\hss}%
\def\getG(#1;#2;#3;#4){$(#1\scolon#2\comma#3\scolon#4)$\hss}%
\centerline{\vbox{\halign{\hss$#$&\quad\hss$#$\hss&\quad\get#&
 \quad\hss\getG#\hss\cr
\omit\hss$p$\hss&$N$&\omit\qquad
 Factors $\mp_\xi\in\FF{p}[t]$ of $\Delta$\hss&\omit\quad\hss$\bG\subset\MG$\hss\cr
\noalign{\vskip2pt \hrule}\cr
  2&\*15&[t^4+t+1, t^4+t^3+1]&(17;1;2;1^2 15^1)\cr
  5&  12&[t^2+2t+4, t^2+3t+4]&(52;0;4;1^4 12^4)\cr
 13&\*12&[t+2, t+6, t+7, t+11]&(14;0;2;1^2 12^1)\cr
 19&  18&[t+2, t+3, t+10, t+13, t+14, t+15]&(40;2;4;1^2 2^1 18^2)\cr
\crcr}}}\egroup
\endtable
\endtheorem

\addendum\label{add.N<=5}
In the settings of \autoref{th.p}, assume in addition that $N\le5$ and
$\CC$ is irreducible.
Then
the pair $(p,N)$
and\mnote{$\rank\to\dim$, statement on $\dim$, edited}
the dimension $r:=\dim_{\bbbk}\AM_\CC(\xi)$
can take one of the following values\rom:
\roster
\item
$(p,N)=(3,4)$ and $r=1$\rom;
\item
$(p,N)=(3,1)$ or $(p,3)$, $p\ne3$, with $r\le2$\rom;
\item
$(p,N)=(5,1)$ or $(p,5)$, $p\ne5$, with $r\le2$\rom;
\item
$(p,N)=(7,1)$ and $r=1$.
\endroster
All four possibilities\mnote{extended} for $(p,N)$ \rom(and all possibilities
for~$r$\rom)
are realized by \emph{genuine} trigonal curves,
and for such curves they are mutually exclusive.
\endaddendum

\addendum\label{add.N>10}
For each pair $(p,N)$ as in \autoref{tab.factors}, at most one of
the factors $\mp_\xi$
listed can appear in the Alexander polynomial of any
given curve.
The six triples $(p,N,\mp_\xi)$ marked with a $^*$ in the table
do appear in the Alexander
polynomials of \emph{genuine} trigonal curves\rom;
the other eight do not.
\endaddendum

\autoref{th.0} is proved in Subsection~\ref{proof.0}.
\autoref{th.p} and \autoref{add.N>10} are proved in
Subsection~\ref{proof.p}, and \autoref{add.N<=5} merely
summarizes the detailed description of the modules
$\AM_\CC/\cp_N(-t)$, $N\le5$, given in
Subsections~\ref{s.N=3}--\ref{s.many}.

In \autoref{tab.factors}, the last column gives a description of
the projection to the modular group $\MG:=\PSL(2,\Z)$ of the
corresponding universal subgroup, see
\autoref{def.universal}. Listed are the index $[\MG:\bG]$,
the numbers~$c_2$, $c_3$ of the conjugacy classes of elements of
order~$2$ and~$3$, respectively, and the set of cusp widths in the
partition notation, see~\cite{Cummins.Pauli}. These data do not
determine the subgroup completely, but drawing large diagrams does
not seem practical here. Note that, in each case marked with
a~$^*$,
the universal subgroup~$\bG'$
corresponding to \emph{genuine} trigonal curves is smaller than
the one listed: one has $[\bG:\bG']=3$. Each time, the skeleton
of~$\bG$, see Subsection~\ref{s.Sk}, has one monovalent
\black-vertex and one monogonal region with the type specification
nontrivial modulo~$6$, see Subsection~\ref{s.type}, and the
skeleton of~$\bG'$ is the triple cyclic covering ramified at these
vertex and region.

\subsection{Ramifications and speculations}
The assumption that the trigonal curve in question
should be irreducible is not very
important. Lifting this requirement would result in a few extra
factors with $N=1$, $2$, or~$4$; they are controlled by
congruence subgroups and thus can easily be enumerated, see
Subsections~\ref{s.N=1,2,4} and \ref{s.cong}.
(The case $N=1$ is known, see~\cite{degt:trigonal}.)

As an extra addendum, mention that, for genuine trigonal curves,
each triple $(p,N,\mp_\xi)$ among
those listed appears in the Alexander polynomial `in a unique
way', in the sense that,
up to Nagata equivalence, each curve~$\CC$ with
$\mp_\xi\divides|\Delta_{\CC,p}$ is
induced from a certain universal curve with this property, see
Subsection~\ref{s.curves} for the definitions.
This statement follows from the uniqueness of the
corresponding universal subgroups (found in the computation)
and \autoref{th.monodromy}.

All four statements apply equally well to plane curves with a
singular point of multiplicity ${\deg}-3$ (as they can be regarded
as generalized trigonal curves in the Hirzebruch
surface~$\Sigma_1$, see Subsection~\ref{s.plane}),
provided that the trigonal model of the curve is not isotrivial. The
relatively simple case of
irreducible isotrivial curves is discussed in
Subsection~\ref{s.isotrivial}; the degree of the Alexander
polynomials of such curves is not universally bounded.

The parabolic case $N=6$ is treated in Section~\ref{S.N=6}; we do
not mention it here as it does not seem to lead to nontrivial
\emph{conventional} Alexander polynomials. (In fact, we mainly
study the so-called \emph{extended} Alexander polynomials,
which depend on the monodromy
group of the curve rather than on its fundamental group only, see
\autoref{def.eAP}
and \autoref{rem.eAP}.)\mnote{extra ref}
The range $7\le N\le10$ remains open. A
few examples are found in \autoref{tab.examples}
in \autoref{example.7--10}. I conjecture that
Tables~\ref{tab.factors} and~\ref{tab.examples} do exhaust all
possibilities with $N\ge7$.
Among other consequences,
this conjecture
would imply that, as an abelian group, $\AM_\CC$ may have
$p$-torsion for finitely many primes~$p$ only; the current list is
$2\le p\le43$ but $p\ne23$, $31$, or~$41$.

Another question
left open for $N>5$ is which pairs,
triples, \etc\. of factors~$\mp_\xi$ can appear simultaneously in
the Alexander polynomial of a particular curve. This problem
reduces to computing the genera of the intersections of the
corresponding universal subgroups, including all their conjugates,
or, equivalently, the genera of
the connected components of
the fibered products of their
skeletons. We postpone this computation until the conjecture
above has been settled.

It is
worth mentioning that none of the groups~$\bG$ listed in
Tables~\ref{tab.factors} and~\ref{tab.examples}
is a congruence subgroup of~$\MG$ (which is easily shown using the
`signatures' listed
and
the tables found in~\cite{Cummins.Pauli}).
This fact refutes my
original expectation that the fundamental group of a
non-isotrivial genuine trigonal
curve might be controlled by congruence subgroups.

\subsection{Idea of the proof}
Modifying the classical Zariski--van Kampen theorem,
see \autoref{th.vanKampen},
one reduces the study of the fundamental group of a (generalized)
trigonal curve~$\CC$ to a question about its monodromy group
$\BM_\CC$, which is a subgroup of the braid group~$\BG3$
(respectively, of its extension \via\ the inner automorphisms of
the free group~$\AA$).
Crucial is the fact that the projection of $\BM_\CC$ to the
modular group~$\MG$ is
a subgroup of genus zero, see~\cite{degt:trigonal} and
\autoref{th.monodromy},
which imposes a very strong restriction to $\BM_\CC$.
The Alexander polynomial is controlled by the reduced Burau
representation, see \cite{Burau,Libgober:Burau} and
Subsection~\ref{s.Burau}, which is
a $\BG3$-action on a certain universal Alexander
$\Lambda$-module $\AM\cong\Lambda\oplus\Lambda$.
Then,
it remains to describe the `Burau congruence subgroups'
$\{\Gb\in\BG3\,|\,\Gb={\id}\bmod\CV\}$, where $\CV\subset\AM$ is a
fixed submodule, and select those that are of genus zero.

Unfortunately, no convenient description of the image of~$\BG3$ in
$\MAT(\Lambda)$ seems to be known, and we choose a more geometric
approach. A subgroup $G\subset\BG3$ is represented by its
\emph{skeleton}~$\Sk$, see Subsection~\ref{s.Sk}, which is a
certain planar (in the case of genus zero) bipartite ribbon graph. Then, in
Section~\ref{S.monodromy}, we derive some local restrictions to the
geometry of~$\Sk$ necessary for the nonvanishing of the Alexander
module. In Section~\ref{S.proof.p}, these local restrictions and
the planarity condition
(Euler's\mnote{corrected} formula $\chi(S)=2$, where $S$ is the minimal
supporting surface of~$\Sk$)
are used to narrow~$N$ down
to the range $N\le26$ (or $N\le21$ if $p=0$). In this finite
range, we use a computer aided analysis to improve the
\latin{a priori} bound on the number of `small' regions of~$\Sk$
and reduce it further to $N\le10$, with the exception of finitely
many triples $(p,N,\mp_\xi)$, $p>0$, see \autoref{N<=10}.\mnote{ref added}
For each exceptional triple, we compute the genus of the
corresponding universal subgroup~$G$ by a straightforward coset
enumeration in the finite group $\GL(2,\bbbk_p)$, thus proving
\autoref{th.p}.

In Section~\ref{S.proof.0}, the case $N\le5$ is reduced to
congruence subgroups of~$\MG$, allowing for an easy classification
of the Alexander modules. Then, for $p=0$, we eliminate the range
$6\le N\le10$ and prove \autoref{th.0}. (For $N=7$ and~$9$, we
have to use {\tt Maple} to show that the corresponding universal
subgroups are of infinite index.)

Sections~\ref{S.BG} and~\ref{S.curves} are preliminary: we
introduce the groups used and necessary technical tools and
explain the relation between trigonal curves and subgroups
of~$\BG3$. Section~\ref{S.N=6} deals with the parabolic case
$N=6$: we discover an infinite series of non-congruence subgroups
of genus zero with nontrivial extended Alexander module.

\remark\label{rem.Burau}
As an interesting by-product of this research,
not quite related to the original problem,
we discover that
`Burau congruence subgroups' described above behave quite
differently from the conventional congruence subgroups of~$\MG$:
there are finitely many subgroups for $N\le5$, infinitely many
finite index subgroups, all of genus zero or one, for $N=6$, and
the subgroups seem to be of infinite index for $N\ge7$ (although
formally the latter claim has only been proved for $N=7$ and~$9$).
Apparently, this is due to the fact that the Burau representation
on $\AM/\cp_N(-t)$ is highly nontransitive for $N\ge7$.
\endremark

\subsection{Acknowledgements}
I am grateful to A.~Libgober for his helpful remarks and
stimulating discussions of the subject.
The\mnote{extended}
final version of the manuscript was prepared during my sabbatical stay at
\emph{l'Instutut des Hautes \'{E}tudes Scientifiques} and
\emph{Max-Planck-Institut f\"{u}r Mathematik};
I would like to extend my gratitude to
these institutions for their support and hospitality.

\section{The braid group\label{S.BG}}

In this section, we introduce the braid group~$\BG3$ and related
objects, the principal purpose being fixing the notation and
terminology.

\subsection{The group $\BG3$}
Let $\AA=\<\Ga_1,\Ga_2,\Ga_3\>$ be the free group on three
generators. The \emph{braid group} $\BG3$ can be defined as the
group of automorphisms $\Gb\:\AA\to\AA$ with the following
properties:
\roster*
\item
each generator~$\Ga_i$ is taken to a conjugate of a generator;
\item
the element $\Gr:=\Ga_1\Ga_2\Ga_3$ remains fixed.
\endroster
Recall, see~\cite{Artin}, that
$\BG3=\<\Gs_1,\Gs_2\,|\,\Gs_1\Gs_2\Gs_1=\Gs_2\Gs_1\Gs_2\>$, the
\emph{Artin generators} $\Gs_1$, $\Gs_2$ acting on~$\AA$ \via
\[*
\Gs_1\:\Ga_1\mapsto\Ga_1\Ga_2\Ga_1\1,\quad
 \Ga_2\mapsto\Ga_1;\qquad
\Gs_2\:\Ga_2\mapsto\Ga_2\Ga_3\Ga_2\1,\quad
 \Ga_3\mapsto\Ga_2.
\]
Note that the set of Artin generators depends on the
basis $\{\Ga_1,\Ga_2,\Ga_3\}$.

In the sequel, we reserve the notation $\AA$ for the free group
supplied with a $\BG3$-action, or, equivalently, with a
distinguished set of bases constituting a whole $\BG3$-orbit. Any
basis in the distinguished orbit is called \emph{geometric}; any
such basis gives rise to a pair of Artin generators of~$\BG3$.
We will also consider the \emph{degree homomorphisms}
\[*
\deg\:\AA\to\Z,\quad\Ga_1,\Ga_2,\Ga_3\mapsto1,\qquad
\bdeg\:\BG3\to\Z,\quad\Gs_1,\Gs_2\mapsto1.
\]
It is straightforward that they do not depend on the choice of
a geometric basis $\{\Ga_1,\Ga_2,\Ga_3\}$ and that for any
$\Ga\in\AA$, $\Gb\in\BG3$ one has $\deg\Gb(\Ga)=\deg\Ga$.

With generalized trigonal curves in mind, see
Subsection~\ref{s.generalized},
introduce also the extended group
$\BG3\cdot\Inn\AA\subset\Aut\AA$, where
$\Inn\AA\cong\AA$ is the subgroup of the inner
automorphisms of~$\AA$. The intersection $\BG3\cap\Inn\AA$ is
the cyclic group generated by $(\Gs_2\Gs_1)^3=\Gr$; hence the
degree map extends to the product \via\
$\bdeg(\Gb\cdot\Ga)=\bdeg\Gb+2\deg\Ga$, where $\Gb\in\BG3$ and
$\Ga\in\Inn\AA\cong\AA$.

The natural action of $\Prod$ on the set of conjugacy classes of
geometric generators defines an epimorphism $\Prod\onto\SG3$. A
subgroup $G\subset\Prod$ is
said to be
\emph{$\SS$-transitive} if this
action, restricted to~$G$, is transitive. Clearly, $G$ is
$\SS$-transitive if and only if its image under the above
epimorphism contains
a cycle of length three.

Given two subgroups~$G$, $H$ of $\BG3$ or $\Prod$
(or any of the quotients
$\Bu3$, $\tMG$, or~$\MG$ considered below), we write $G\sim H$ if
$G$ is conjugate to~$H$ and $G\prec H$ if $G$ is
\emph{subconjugate} to~$H$, \ie, if $G$ is conjugate to a subgroup
of~$H$.

\subsection{The Burau representation\label{s.Burau}}
Denote by~$\AM$ the abelianization of the kernel $\Ker\deg$, and
let $[h]\in\AM$ be the class of an element $h\in\Ker\deg$. An
element $\Ga\in\AA$ of degree one defines a homomorphism
$t\:\AM\to\AM$, $[h]\mapsto[\Ga h\Ga\1]$, which does not depend
on~$\Ga$. Thus, $\AM$ turns into a module over the ring
$\Lambda:=\Z[t,t\1]$ of Laurent polynomials. An easy computation
shows that $\AM=\Lambda\be_1\oplus\Lambda\be_2$, where
$\be_1=[\Ga_2\Ga_1\1]$, $\be_2=[\Ga_3\Ga_2\1]$ in some geometric
basis $\{\Ga_1,\Ga_2,\Ga_3\}$.

Since the $\BG3$-action on~$\AA$ preserves the degree, it restricts to
a certain action on~$\AM$, which is
called the \emph{\rom(reduced\rom) Burau
representation}, see~\cite{Burau}.
This representation is faithful; for this reason we identify an
element $\Gb\in\BG3$ and the matrix in
$\MAT(\Lambda)$
representing it.
The Artin generators~$\Gs_1$, $\Gs_2$
corresponding to the chosen
geometric basis $\{\Ga_1,\Ga_2,\Ga_3\}$ (the one used to
define~$\be_1$, $\be_2$)
act \via
\[*
\Gs_1=\mat[-t,1,0,1],\quad
\Gs_2=\mat[1,0,t,-t],
\]
and the powers of these matrices are given by
\[
\Gs_1^m=\mat[(-t)^m,\tf_m(-t),0,1],\quad
\Gs_2^m=\mat[1,0,t\tf_m(-t),(-t)^m],
\label{eq.powers}
\]
where $\tf_m(t):=(t^m-1)/(t-1)$.
For future references, observe that, for any $r\in\Z$,
one has
\[
(t+1)t^{r}\tf_m(-t)+t^{r}(-t)^m=t^{r}.
\label{eq.identity}
\]
The following two matrices are also used in the sequel:
\[*
\Gs_2\Gs_1=\mat[-t,1,-t^2,0],\quad
\Gs_2\Gs_1\Gs_2=\mat[0,-t,-t^2,0].
\]

The Burau representation extends to the product
$\Prod$. Clearly, the map
$\Inn\AA=\AA\to\MAT(\Lambda)$ is given
by $\Ga\mapsto t^{\deg\Ga}\id$. The image of
$\Prod$ in the group $\GL(2,\Lambda)$
is denoted by $\Bu3$; it is the central product
$\BG3\odot\Z$, obtained by identifying the center $Z(\BG3)$
and the subgroup $3\Z\subset\Z$ (both
subgroups being generated by $t^3\id$).
The center $Z(\Bu3)$ is the cyclic subgroup formed by all scalar
matrices $t^r\id$. The degree map $\bdeg$ descends to~$\Bu3$ and
coincides, essentially, with the determinant: one has
$\det\Gb=(-t)^{\bdeg\Gb}$ for any $\Gb\in\Bu3$.

Given two submodules $\CU,\CV\subset\AM$, we say that $\CU$ is
\emph{conjugate} to~$\CV$, $\CU\sim\CV$, if $\CV=\Gb(\CU)$ for
some $\Gb\in\BG3$, and $\CU$ is \emph{subconjugate} to~$\CV$,
$\CU\prec\CV$, if $\CU$ is conjugate to a submodule of~$\CV$.
Clearly, in this definition $\BG3$ can be replaced with~$\Bu3$.

For an ideal $I\subset\Lambda$, we will use the notation
$\CU\sim\CV\bmod I$ and $\CU\prec\CV\mod I$ meaning the images of
the modules in $\AM/I$. If $I=\Lambda f$, $f\in\Lambda$,
is a principal ideal, we abbreviate $\bmod\,\Lambda f$ to
$\bmod\,f$.

\subsection{The modular representation}\label{s.modular}
Specializing all matrices at $t=-1$, one obtains homomorphisms
$\BG3,\Bu3\to\tMG:=\SL(2,\Z)$, which give rise to the \emph{modular
representation}
\[*
\pr_\MG\:\BG3,\Bu3\to\MG:=\PSL(2,\Z)=\tMG/\!\pm\id.
\]
Usually, we abbreviate $\pr_\MG G=\bG$ and $\pr_\MG\Gb=\Bb$
for a subgroup $G\subset\Bu3$ or an element $\Gb\in\Bu3$.

Recall that the modular group~$\MG$ is generated by two elements
$\X$, $\Y$ subject to the relations $\X^3=\Y^2=1$. One can take
$\X=(\Bs_2\Bs_1)\1$ and $\Y=\Bs_2\Bs_1^2$;
then $\Bs_1=\X\Y$ and $\Bs_2=\X^2\Y\X\1$.

A subgroup of~$\MG$ is called a \emph{congruence subgroup} of
level $l\divides|n$ if it contains the \emph{principal congruence
subgroup} $\MG(n)=\{g\in\MG\,|\,g={\id}\bmod n\}$. We make use
of the list of congruence subgroups found
in~\cite{Cummins.Pauli};
when referring to such subgroups, we use the notation
of~\cite{Cummins.Pauli} and, whenever available, the alternative
conventional notation.

The degree homomorphisms $\bdeg\:\BG3\to\Z$ and $\bdeg\:\Bu3\to\Z$
descend to well defined homomorphisms
$\bdeg\:\MG\to\CG6$ and ${\bdeg}\bmod2\:\MG\to\CG2$, respectively.
Thus,
one has
$\BG3=\MG\times_{\CG6}\Z$ and $\Bu3=\MG\times_{\CG2}\Z$.

\definition\label{def.depth}
The \emph{depth} $\depth G$ of a subgroup $G\subset\Bu3$ is the
degree of the positive generator of the intersection
$G\cap\Ker\pr_\MG$, or zero if this intersection is trivial.
One has $\depth G=0\bmod2$ and $\depth G=0\bmod6$ if
$G\subset\BG3$.
\enddefinition

Consider a subgroup $G\subset\Bu3$, let $2d=\depth G$, and let
$G_d$ be the image of~$G$ under the projection
$\pr_d:={\pr_\MG}\times({\bdeg}\bmod2d)\:\Bu3\to\MG\times\CG{2d}$.
(We let $\CG0=\Z$.) Then $G=\pr_d\1G_d$ and $G_d$ projects
isomorphically onto $\bG$; in other words, $G_d$
is the graph of a certain homomorphism
$\Gf\:\bG\to\CG{2d}$.
This construction is summarized by the following definition and
proposition.

\definition\label{def.slope}
The homomorphism $\Gf\:\bG\to\CG{2d}$ as above is called the
\emph{slope} of a subgroup $G\subset\Bu3$.
\enddefinition


\proposition\label{G<->slope}
There is a one-to-one correspondence between the set of subgroups
$G\subset\Bu3$ and
the set of pairs $(\bG,\Gf)$, where $\bG\subset\MG$
is a subgroup and $\Gf$ is a homomorphism $\bG\to\CG{2d}$ with
the property $\Gf={\bdeg}\bmod2$. One has $G\subset\BG3$ if and
only if $d=0\bmod3$ and $\Gf={\bdeg}\bmod6$.
\qed
\endproposition

Each subgroup $\bG\subset\MG$ admits three canonical slopes,
namely, the restrictions to~$\bG$ of
the homomorphisms $\pm\bdeg\:\MG\to\CG6$ and
${\bdeg}\bmod2\:\MG\to\CG2$. We denote the corresponding
subgroups of~$\Bu3$
by $(\bG)^\pm$ and $(\bG)\two$, respectively.
The subgroups $(\bG)\two=\pr_\MG\1\bG$ and
$(\bG)\plus=(\bG)\two\cap\BG3$ are merely the full preimages
of~$\bG$ under
$\pr_\MG\:\Bu3\to\MG$ and $\pr_\MG\:\BG3\to\MG$, respectively.

\subsection{Skeletons}\label{s.Sk}
In this subsection, we outline the relation between
subgroups of~$\MG$ and certain bipartite ribbon graphs, called
skeletons.
This\mnote{refs added}
and other very similar constructions have been studied, \eg,
in~\cite{Birch,Bogomolov,Kulkarni}. In the exposition below we
follow recent paper~\cite{degt:monodromy},
where all proofs and further details can be found.

Recall\mnote{new paragraph: geometric insight added}
that a \emph{bipartite graph} is a graph whose vertices are divided into two
kinds, \black-- and \white--, so that the two ends of each edge are of the
opposite kinds. A \emph{ribbon graph} is a graph equipped with a
distinguished cyclic order (\iq. transitive $\Z$-action) on the star of each
vertex. Any graph embedded into an oriented surface~$S$ is a ribbon graph,
with the cyclic order induced from the orientation of~$S$. Conversely, any
finite ribbon graph defines a unique, up to homeomorphism, closed oriented
surface~$S$ into which it is embedded: the star of each vertex is embedded
into a small oriented disk (it is this step where the cyclic order is used),
these disks are connected by oriented ribbons along edges
producing a tubular neighborhood of the graph, and finally each
boundary component of the resulting compact surface is patched with a disk.
(Intuitively, the boundary components patched at the last step are
the \emph{regions} defined combinatorially in
Subsection~\ref{ss.regions} below.)
The surface~$S$ thus constructed is called the \emph{minimal supporting
surface} of the ribbon graph.


In the rest of this section,\mnote{edited}
we redefine a certain class of bipartite ribbon
graphs in purely combinatorial terms, relating them to the modular group.
In spite of this combinatorial approach, we will
freely use the topological language applicable to the geometric
realizations of the graphs.

\paragraph\label{ss.skeleton}
Given a subgroup $G\subset\MG$, its \emph{skeleton} $\Sk=\Sk_G$ is
the bipartite ribbon graph, possibly infinite, defined as follows:
the set of edges of~$\Sk$ is the $\MG$-set
$\MG/G$, its \black-- and \white-vertices are the orbits of~$\X$
and~$\Y$, respectively, and the cyclic order (ribbon graph
structure) at a trivalent \black-vertex is given by~$\X\1$. (All
other vertices are at most bivalent and cyclic order is
irrelevant.) The skeleton $\Sk_G$ is equipped with a distinguished
edge, namely, the coset $G/G$.

By definition, $\Sk$ is a connected bipartite graph with
the following properties:
\roster*
\item
the valency of each \black-vertex equals $1$ or~$3$
(a divisor of $\ord\X=3$),\mnote{explained} and
\item
the valency of each \white-vertex equals~$1$ or~$2$
(a divisor of $\ord\Y=2$).\mnote{explained}
\endroster
Conversely, the set of edges of any connected
bipartite ribbon graph~$\Sk$
satisfying the valency restriction above
admits a natural structure of a transitive $\MG$-set
(the action of~$\X\1$ and~$\Y$ following the cyclic order at the \black-- and
\white-vertices, respectively),\mnote{explained}
and the
original subgroup~$G$ can be recovered, up to conjugation, as the
stabilizer $\Stab(e)$
of any edge~$e$ of~$\Sk$.

\convention
In the figures,
we omit
bivalent \white-vertices, assuming that such a vertex is to be
inserted at the center of each edge connecting two
\black-vertices. With an abuse of the language, we will speak
about
\emph{adjacent} \black-vertices, meaning that they are
connected by
a pair of edges with a common bivalent \white-vertex.

As usual,
skeletons of genus zero (see Subsection~\ref{ss.regions} below)\mnote{ref added}
are drawn in the disk, assuming the
blackboard thickening for the ribbon graph structure.
The boundary of the disk (the dotted grey
circle in the figures) represents a single point in the
sphere~$S^2$.
\endconvention

\paragraph\label{ss.path}
Topologically, it is convenient to regard~$\Sk$ as an orbifold,
assigning to each monovalent \black-- or \white-vertex
ramification index~$3$ or~$2$, respectively. Then there is a
canonical isomorphism
\[*
G=\Stab(e)=\piorb(\Sk,e),
\]
where the
basepoint for the fundamental group is chosen inside an edge~$e$.
In fact, homotopy classes of
paths in $\Sk$
(taking into account the
orbifold structure) can be identified with pairs $(e_0,g)$, where
the \emph{starting point}~$e_0$ is an edge and $g\in\MG$; the
\emph{ending point} of such a path is then $e_1:=g\1e_0$.
Intuitively,\mnote{geometric insight explained}
one starts at~$e_0$ and constructs a path edge by edge, choosing at each
steps between one of the four possible directions:
turning about the \black-- or \white-end
of the last edge in the positive or negative direction (with respect to the
distinguished cyclic order); these directions are encoded by
the letters $\Y\1=\Y$ or
$\X^{\minusplus1}$ in the word representing~$g$.

A path $(e,g)$, $g\in\MG$,
is a loop if and only if $e=g\1e$, \ie, $g\in\Stab(e)$;
hence the isomorphism above.

\paragraph\label{ss.regions}
A \emph{region} of a skeleton~$\Sk$ is an orbit of~$\X\Y$. The
cardinality of a region~$R$ is called its \emph{width} $\rdeg R$.
(In the arithmetical theory, instead of regions one speaks about
cusps and cusp widths; this, and the fact that the term `degree'
is way too overused,
explains the terminology.) A region~$R$
of
width~$n$ is also referred to as an \emph{$n$-gon} or
\emph{$n$-gonal region}, `corners' being the \black-vertices in
the boundary of~$R$. If $\Sk$ is finite, then, patching each
region with an oriented disk, one obtains a minimal compact
oriented surface~$S$ supporting~$\Sk$. Its genus is called the
\emph{genus} of~$\Sk$ and of the subgroup $G\subset\MG$
corresponding to~$\Sk$. (This definition is equivalent to the
conventional one, see~\cite{degt:monodromy}.)
Using the projection $\pr_\MG$, we extend
the notions of skeleton, genus, \etc. to subgroups of~$\Bu3$.

A \emph{marking} at a trivalent \black-vertex~$v$ is a choice of
an edge~$e$ adjacent to~$v$.
The region (orbit) containing an edge~$e$ is denoted by $\reg(e)$.
Thus, the three
regions
adjacent to a marked vertex $(v,e)$
are $\reg(e)$, $\reg(\X e)$, and
$\reg(\X^2e)$.
By default, given a region~$R$,
a marking~$e$ at each vertex~$v$ in $\dR$
is chosen so that $R=\reg(e)$.
Note that a vertex
may appear in~$\dR$ more then once; in this case each occurrence
gets its own marking.

\paragraph\label{ss.coverings}
An inclusion $G'\subset G$ of two subgroups gives rise to a
$\MG$-map $\Sk'\to\Sk$ of their skeletons, which is
a covering with respect to the orbifold structure defined
in Subsection~\ref{ss.path}.
It extends to an essentially unique (ramified)
covering $S'\to S$ of
the minimal surfaces, see Subsection~\ref{ss.regions}.
The covering $\Sk'\to\Sk$ is called
\emph{\rom(un-\rom)ramified} if so is $S'\to S$. In other words,
the covering is unramified if and only if the pull-back of each
monovalent vertex of~$\Sk$ consists of monovalent vertices only
and the
pull-back of each region~$R$ of~$\Sk$ consists of regions of
the same width $\rdeg R$.

\paragraph\label{ss.basis}
In the definition of the skeleton~$\Sk$ of a subgroup
$G$, we use a
distinguished pair $\X$, $\Y$ of generators of~$\MG$, hence a
distinguished pair $\Gs_1$, $\Gs_2$ of Artin generators of~$\BG3$,
hence a distinguished geometric basis $\{\Ga_1,\Ga_2,\Ga_3\}$
of~$\AA$; the latter is defined up to the action of the center
$Z(\BG3)$, \ie, up to conjugation by~$\Gr$.

One has $G=\piorb(\Sk,e)$, where $e=G/G$ is the distinguished
edge of~$\Sk$, see Subsection~\ref{ss.path}.
If $e'$ is another edge, we fix a path $\Gg=(e,g)$ from~$e$
to~$e'$ and identify $\piorb(\Sk,e')$ with~$G$ \via\ the
translation isomorphism $\Gd\mapsto\Gg\Gd\Gg\1$, \ie, \via\ the
conjugation by~$g$.
Alternatively, one can lift~$g$ to an element $\tilde g\in\BG3$
and consider the new geometric basis
$\{\Ga_1',\Ga_2',\Ga_3'\}$, $\Ga_i'=\tilde g(\Ga_i)$, for~$\AA$.
In this sense, assuming~$\Gg$ fixed, we will speak about a
\emph{canonical basis} over~$e'$.

\subsection{Type specification}\label{s.type}
If $G\subset\Bu3$ is a subgroup of genus zero, its slope can be
described in terms of its skeleton~$\Sk$. In view
of Subsection~\ref{ss.path}, the projection $\bG\subset\MG$ has
a presentation of the form
\[
\textstyle
\bigl<\Gb_R,\Gg_v\bigm|
 (\Gg^\BLACK_v)^3=(\Gg^\WHITE_v)^2=1,\
 \prod\Gb_R\prod\Gg_v=1\bigr>,
\label{eq.G}
\]
where the indices~$R$ and~$v$ run,
respectively, over all regions and monovalent vertices of~$\Sk$
and the superscript indicates the type of the vertex. (The product
in the last relation is in a certain order depending on the
choice of the basis. In fact, $\{\Gb_R,\Gg_v\}$ is merely a
geometric basis for the fundamental group of a punctured sphere,
\cf. \autoref{def.basis} below.)
Furthermore, each generator $\Gb_R$ is
conjugate to $\Bs_1^{\rdeg R}$, and each
generator~$\Gg_v$ is
conjugate to $\X\1=\Bs_2\Bs_1$ or~$\Y=\Bs_2\Bs_1^2$,
depending on whether $v$ is a \black-- or \white-vertex,
respectively.

\definition\label{def.type}
The \emph{type specification} of a subgroup $G\subset\Bu3$ of genus
zero is the $\CG{\depth G}$-valued
function~$\type$ defined on the set of all regions and monovalent
vertices of the skeleton $\Sk_G$; each region or monovalent vertex
is sent to the degree of (any) lift to~$G$ of the corresponding
generator in~\eqref{eq.G} or, equivalently, to the value of the
slope of~$G$ on the corresponding generator.
\enddefinition

\proposition\label{type}
Let $d=6$ if $G\subset\BG3$ and $d=2$ otherwise. Then one has\rom:
\roster
\item\label{tp.dp}
$\depth G=0\bmod d$\rom;
\item\label{tp.R}
$\type(R)=\rdeg R\bmod d$ for any region~$R$\rom;
\item\label{tp.black}
$\type(\BLACK)=2\bmod d$ and $3\type(\BLACK)=0$\rom;
\item\label{tp.white}
$\type(\WHITE)=3\bmod d$ and $2\type(\WHITE)=0$\rom;
\item\label{tp.sum}
the sum of all values of~$\type$ equals zero.
\endroster
Any pair $(\depth,\type)$ satisfying~\loccit{tp.dp}--\loccit{tp.sum}
above
defines a unique slope\rom;
such a pair results
in a subgroup $G\subset\BG3$ if and only if it
satisfies~\loccit{tp.dp}--\loccit{tp.white} with $d=6$.
\endproposition

\proof
The $({\bmod}\,d)$-congruences
in~\loccit{tp.dp}--\loccit{tp.white}
follow from the properties of
slopes, see \autoref{G<->slope}, and the other relations
in~\loccit{tp.black}--\loccit{tp.sum} are the abelian versions of
the relations in~\eqref{eq.G}. The type specification
determines the slope of~$G$ as it assigns a value to each generator
in~\eqref{eq.G}.
\endproof

Given an integer~$m$, a type specification is said to be
\emph{trivial modulo~$m$} if it satisfies the congruences
in~\iref{type}{tp.dp}--\ditto{tp.white} with $d=m$. Thus,
\autoref{type} states that any type specification is
trivial modulo~$2$ and that a subgroup~$G$ is in~$\BG3$
if and only if its type
specification is trivial modulo~$6$.

\convention
In the drawings,
we indicate the type specification (inside
a region or next to a vertex) only when it is not trivial
modulo~$0$.
\endconvention

\subsection{The Alexander module}\label{s.kernels}
For a subgroup $G\subset\Prod$, let
\[*
\textstyle
\IM_G=\sum_{\Gb\in G}\Im(\Gb-\id)\subset\AM,\qquad
\tIM_G=\sum_{\Gb\in G,\ \Ga\in\AA}\Lambda[\Gb(\Ga)\cdot\Ga\1]
 \subset\AM.
\]

\definition\label{def.AP}
The \emph{Alexander module} of a subgroup $G\subset\Prod$
is the $\Lambda$-module $\AM_G:=\AM/\tIM_G$.
If the product $\AM_G\otimes\Bbbk_p$ is a torsion
$(\Lambda\otimes\Bbbk_p)$-module, its
order $\Delta_{G,p}\in\Lambda\otimes\Bbbk_p$ is called the
\emph{$({\bmod}\,p)$-Alexander polynomial} of~$G$.
We usually abbreviate $\Delta_{G,0}=\Delta_G$.
\enddefinition

\definition\label{def.eAP}
The \emph{extended Alexander module} of a subgroup $G\subset\Prod$
is the $\Lambda$-module $\bAM_G:=\AM/\IM_G$; the
\emph{extended Alexander polynomial}
$\bDelta_{G,p}\in\Lambda\otimes\Bbbk_p$
(whenever defined) is the order
of the $(\Lambda\otimes\Bbbk_p)$-module $\bAM_G\otimes\Bbbk_p$.
\enddefinition

Clearly, the Alexander polynomial $\Delta_{G,p}$
and its extended counterpart~$\bDelta_{G,p}$
can be computed using any field~$\bbbk$ of characteristic~$p$,
and the Alexander polynomial can be
interpreted as the characteristic polynomial of the operator~$t$
acting on the finite dimensional $\bbbk$-vector space
$\AM_G\otimes\bbbk$ (respectively, $\bAM_G\otimes\bbbk$).

\remark\label{rem.eAP}
Assume that $G=\BM_C$\mnote{remark moved here, entitled, and extended}
is the monodromy group of a trigonal curve, see
Subsections~\ref{s.sections} and~\ref{s.generalized}
below.
Then, the conventional Alexander module $\AM_G$ is the Alexander module
of~$C$; it depends on the fundamental group of~$C$ only, see
Subsection~\ref{s.pi1}.
On the contrary,
the submodule $\IM_G\subset\AM$ depends only on the
image of~$G$ in~$\Bu3$; thus, it is easier to compute. Furthermore,
$\IM_G$, $\bAM_G$, and the
extended
Alexander polynomials can be \emph{defined} for subgroups
$G$ of~$\Bu3$ rather than those of the more complicated group $\Prod$.
There is a canonical epimorphism $\bAM_G\onto\AM_G$, \cf. \autoref{lem.V},
and the conventional Alexander polynomials divide their extended
counterparts (whenever defined).
For this reason, and since we are mainly interested
in an upper bound on the Alexander polynomial, we will usually
deal with the extended versions.
\autoref{cor.V} and \autoref{cor.V=V} below
show that, for subgroups of~$\BG3$ (\latin{i.q\.} genuine trigonal curves),
the two submodules $\IM_G,\tIM_G\subset\AM$
usually coincide.
\endremark

\definition\label{def.universal}
Given a submodule $\CV\subset\AM$, the set
\[*
G_{\CV}=\bigl\{\Gb\in\Bu3\bigm|\Im(\Gb-\id)\subset\CV\bigr\}
\]
is a subgroup of~$\Bu3$, \cf.~\cite{degt:trigonal}; it is called
the \emph{universal subgroup} corresponding to~$\CV$.
\enddefinition

Definitions~\ref{def.AP}, \ref{def.eAP}, and~\ref{def.universal}
have a geometric meaning for subgroups of genus zero, see
Subsection~\ref{s.pi1} below. In general, it is not quite clear
how the Alexander modules and,
especially, universal subgroups should be defined,
see \autoref{rem.meaning}.

Next two statements are straightforward.

\lemma\label{lem.universal}
For subgroups $G,H\subset\Bu3$ and submodules $\CU,\CV\subset\AM$,
one has
\roster
\item
if $G\prec H$, then $\IM_G\prec\IM_H$\rom;
\item
if $\CU\prec\CV$, then $G_{\CU}\prec G_{\CV}$\rom;
\item
$\IM_G\prec\CU$ if and only if $G\prec G_{\CU}$.
\qed
\endroster
\par\removelastskip
\endlemma

\lemma\label{lem.V}
One has\rom:
\roster
\item\label{V.1}
$\IM_G\subset\tIM_G$,
\item\label{V.2}
$[\Gb(\Ga h)\cdot(\Ga h)\1]=[\Gb(\Ga)\cdot\Ga\1]+t^{\deg\Ga}(\Gb-\id)[h]$
for any $h\in\Ker\deg$,
\item\label{V.3}
$[\Gb(\Ga^n)\cdot\Ga^{-n}]=\tf_n(t^{\deg\Ga})[\Gb(\Ga)\cdot\Ga\1]$
for any $n\in\Z$,
\endroster
where $\Gb\in\BG3$ and $\Ga\in\AA$. As a consequence,
\roster[4]
\item\label{V.4}
$\tIM_G=\IM_G+\sum_{\Gb\in G}\Lambda[\Gb(\Ga_i)\cdot\Ga_i\1]$
for any geometric generator $\Ga_i\in\AA$.
\qed
\endroster
\par\removelastskip
\endlemma

\lemma\label{cor.V}
For a subgroup $G\in\BG3$,
one has $(t^2+t+1)\tIM_G\subset\IM_G$.
\endlemma

\proof
Since $\deg\Ga_1^3=3=\deg\Gr$, for any braid $\Gb\in\BG3$ one has
\[*
(t^2+t+1)[\Gb(\Ga_1)\cdot\Ga_1\1]=[\Gb(\Ga_1^3)\cdot\Ga_1^{-3}]=
[\Gb(\Gr)\cdot\Gr\1]\bmod\IM_G,
\]
see \autoref{lem.V}\loccit{V.3} and~\ditto{V.2}. Since $\Gr$ is
$\BG3$-invariant, this expression is $0\bmod\IM_G$, and the
statement follows from \autoref{lem.V}\loccit{V.4}.
\endproof

\corollary\label{cor.V=V}
For any subgroup $G\subset\BG3$,
field~$\bbbk$, and
polynomial
$f\in\Lambda\otimes\bbbk$ prime to $t^2+t+1$, the images of
$\IM_G$ and $\tIM_G$ in $(\AM\otimes\bbbk)/f$ coincide.
\qed
\endcorollary

\subsection{Specializations}\label{s.specialization}
Recall that we denote $\Bbbk_0=\Q$ and $\Bbbk_p=\FF{p}$ for $p$
prime. If $p$ is understood, we drop the index.
The notation $\mp_\xi\in\Bbbk[t]$ stands for the minimal
polynomial of an element $\xi\ne0$ of an algebraic extension
$\bbbk\supset\Bbbk$.

\definition
The multiplicative order of an element $\xi\in\bbbk^*$
is denoted by $\ord\xi$. (If $\xi$ is not a root of unity, we let
$\ord\xi=\infty$.)
For $N\in\Z_+$ not divisible
by~$p$ (where $p$ is a prime or zero),
introduce $\ee_p(N)$ as follows:
$\ee_2(N)=N$ and
\[*
\ee_p(N)=\begin{cases}
 2N,&\text{if $N=1\bmod2$},\\
 \frac12N,&\text{if $N=2\bmod4$},\\
 N,&\text{if $N=0\bmod4$}
\end{cases}
\]
for $p\ne2$.
Then $\ee_p(\ord\xi)=\ord(-\xi)$
and $\ee_p$ is an involution:
$\ee_p(\ee_p(N))=N$.
\enddefinition

Given~$\xi$ as above, we
define\mnote{`define'; extended}
the \emph{specializations} of~$\Lambda$ and~$\AM$ at~$\xi$ to be
$\Lambda(\xi)=(\Lambda\otimes\Bbbk)/\mp_\xi$ and
$\AM(\xi)=(\AM\otimes\Bbbk)/\mp_\xi$,
respectively.
(The specializations of other relevant modules are defined below on a
case-by-case basis.)
Usually we assume that
$\bbbk=\Lambda(\xi)$; then $\AM(\xi)$ is a $\bbbk$-vector
space of dimension~$2$.

For a subgroup
$G\subset\Prod$, define the \emph{specializations}
$\IM_G(\xi)\subset\tIM_G(\xi)\subset\AM(\xi)$
as the images of, respectively,
$\IM_G\otimes\Bbbk$ and~$\tIM_G\otimes\Bbbk$
in $\AM(\xi)$.
(In general, the maps
$\IM_G\otimes\Bbbk\to\AM\otimes\Bbbk$ are \emph{not} monomorphisms.)
As above, these images can be regarded as
$\bbbk$-vector subspaces. If $G\subset\BG3$ and
$\xi^2+\xi+1\ne0$, the two subspaces coincide, see
\autoref{cor.V=V}.
We denote $\bAM_G(\xi)=\AM(\xi)/\IM_G(\xi)$ and
$\AM_G(\xi)=\AM(\xi)/\tIM_G(\xi)$. The barred versions
of all objects can as well be
defined for a subgroup $G\subset\Bu3$.

We extend the notion of (sub-)conjugacy, see
Subsection~\ref{s.Burau},
and the notation $\sim$
and~$\prec$ to submodules of $\AM\otimes\Bbbk$ and $\AM(\xi)$. The
concept of universal subgroup, see \autoref{def.universal},
can also be extended to submodules of $\AM\otimes\Bbbk$ and
$\AM(\xi)$, and an analog of \autoref{lem.universal} holds
literally.

\section{Trigonal curves\label{S.curves}}

In this section, we introduce (generalized) trigonal curves and
their monodromy groups.
Proofs are mostly omitted; for all details,
see~\cite{degt:trigonal} and references therein.

\subsection{Trigonal curves in Hirzebruch surfaces}\label{s.curves}
A \emph{Hirzebruch surface}~$\Sigma_d$
is a geometrically ruled rational surface with an
\emph{exceptional section}~$\EE$ of self-intersection $-d\le0$.
The \emph{fibers} of~$\Sigma_d$ are the fibers of the ruling
$\Sigma_d\to\Cp1$.
To avoid excessive notation,
we identify fibers and their images in the base~$\Cp1$.
The semigroup of classes of effective divisors on~$\Sigma_d$ is
freely generated by the classes $\ls|E|$ and~$\ls|F|$, where $F$
is any fiber.

A \emph{generalized trigonal curve} is a reduced curve
$\CC\subset\Sigma_d$
intersecting each fiber at three points,
counted with multiplicities; in other words,
$\CC\in\ls|3\EE+3dF|$.
A \emph{\rom(genuine\rom) trigonal curve} is a generalized
trigonal curve
disjoint
from the exceptional section $\EE\subset\Sigma_d$.
A \emph{singular fiber} of a generalized trigonal curve
$\CC\subset\Sigma_d$
is a fiber $F$ of~$\Sigma_d$
intersecting $\CC\cup E$ geometrically at fewer than
four points, \ie, such that either $\CC$ is tangent to~$F$ or
the union $\CC\cup E$
has a singular point in~$F$.

We emphasize that,\mnote{a general remark for \autoref{th.monodromy} and
similar places}
from our point of view, a trigonal curve is
always
a curve
embedded in a certain way to a certain Hirzebruch surface; the latter is
assumed even if not mentioned explicitly. In particular, all (\hbox{iso-,}
\hbox{auto-,}
\etc.) morphisms of trigonal curves are supposed to extend to their
respective surfaces.

The \emph{\rom(functional\rom) $j$-invariant}
$j_\CC\:\Cp1\to\Cp1$ of
a trigonal curve $\CC\subset\Sigma_d$ is the analytic
continuation of the function sending a
nonsingular
fiber~$F$
to the $j$-invariant (divided
by~$12^3$)
of the elliptic
curve covering~$F$
and ramified at
$F\cap(\CC\cup\EE)$.
In\mnote{coordinate description added}
appropriate affine coordinates $(x,y)$ in $\Sigma_d$
(such that $E=\{y=\infty\}$)
the curve~$\CC$ can
be given by its \emph{Weierstra{\ss} equation}
\[*
y^3+3p(x)y+2q(x)=0.
\]
Then
\[*
j_\CC(x)=\frac{p^3}\Delta,\quad\text{where}\quad
 \Delta(x)=p^3+q^2.
\]
The curve~$\CC$
is called \emph{isotrivial} if $j_\CC=\const$.
A\mnote{extra statement}
non-isotrivial trigonal curve~$\CC$ is determined by its $j$-invariant up to
Nagata equivalence, see below.

A positive (negative) \emph{Nagata transformation} is the
birational transformation $\Sigma_d\dashrightarrow\Sigma_{d\pm1}$
consisting in blowing up a point~$P$ on (respectively, not on) the
exceptional section~$\EE$ and blowing down the proper transform of the
fiber through $P$. An \emph{$m$-fold} Nagata transformation is a
sequence of $m$ Nagata transformations \emph{of the same sign}
over the same point of the base.
Two trigonal curves~$\CC$, $\CC'$ are called \emph{$m$-Nagata
equivalent} if $\CC'$ is the proper transform of~$\CC$
under a sequence of $m$-fold Nagata transformations. The
special case $m=1$ is referred to as just \emph{Nagata equivalence}.

Each generalized trigonal curve~$\CC$ is Nagata equivalent to a genuine
one, which is unique up to Nagata equivalence. It is called a
\emph{trigonal model} of~$\CC$.

Given a nonconstant holomorphic map $\tj\:\Cp1\to\Cp1$,
the ruled surface
$\Sigma':=\tj^*\Sigma_d$ is also a Hirzebruch surface;
it is isomorphic to $\Sigma_{d\cdot\deg\tj}$. Given a
trigonal curve $\CC\subset\Sigma_d$, its divisorial pull-back
$\CC':=\tj^*\CC\subset\Sigma'$ is also a trigonal curve; it is
said to be \emph{induced} from~$\CC$ by~$\tj$.

\subsection{Braid monodromy}\label{s.sections}
Introduced\mnote{intro written}
in this subsection are the necessary prerequisites for the
classical Zariski--van Kampen theorem: we define the notion of proper section
and, using such a section, construct the braid monodromy of a curve.
The construction applies literally to any curve disjoint from the exceptional
section;
in the case of a \emph{trigonal} curve $\CC\subset\Sigma_d$,
it turns out that the monodromy group
captures quite a few essential geometric properties of~$\CC$, see
\autoref{th.monodromy} for the precise statement.

Fix a Hirzebruch surface~$\Sigma_d$.
For a fiber~$F$ of~$\Sigma_d$, the complement
$F^\circ:=F\sminus\EE$ is an affine space over~$\C$. Hence, one
can speak about the convex hull of a subset of~$F^\circ$. For a
subset $S\subset\Sigma_d\sminus\EE$, denote by $\Conv_F S$
the convex hull of $S\cap F\scirc$ in~$F\scirc$ and
let $\Conv S=\bigcup_F\Conv_FS$.

Fix a genuine trigonal curve $\CC\subset\Sigma_d$.
The term `section' stands for a continuous section of (a
restriction of) the fibration $p\:\Sigma_d\to\Cp1$.
Let $\disk\subset\Cp1$ be a closed topological disk.
(In what follows, we take for~$\disk$ the complement of a
small regular neighborhood of a nonsingular fiber $F_0\in\Cp1$.)
A section $s\:\disk\to\Sigma_k$ of~$p$ is called \emph{proper}
if its
image is disjoint from both~$\EE$ and $\Conv\CC$.
As a simple consequence of the obstruction theory,\mnote{explained}
any disk $\disk\subset\Cp1$ admits a proper section
$s\:\disk\to\Sigma_k$, unique up to homotopy in the class
of proper sections.

Fix a disk $\disk\subset\Cp1$ and
let $F_1,\ldots,F_r\in\disk$ be all singular and, possibly, some
nonsingular fibers of~$\CC$ that belong to~$\disk$.
Assume that
all these fibers are in the interior of~$\disk$.
Let
$\disk\scirc=\disk\sminus\{F_1,\ldots,F_r\}$ and
fix a \emph{reference fiber}
$F\in\disk\scirc$. Then, given a proper
section~$s$, one can define the group
$\pi_F:=\pi_1(F\scirc\sminus\CC,s(F))$ and the
\emph{braid monodromy}, which is the anti-homomorphism
$\bm\:\pi_1(\disk\scirc,F)\to\Aut\pi_F$ sending a loop~$\Gg$ to
the automorphism obtained by dragging~$F$ along~$\Gg$ and keeping
the reference point in~$s$.

\definition\label{def.basis}
Let~$D$ be an oriented punctured disk, and let $b\in\partial D$. A
\emph{geometric basis} in~$D$ is a basis $\{\Gg_1,\ldots,\Gg_r\}$
for the free group $\pi_1(D,b)$ formed by the classes of
positively oriented lassoes about the punctures,
pairwise disjoint except
at the common reference point~$b$ and such that
$\Gg_1\ldots\Gg_r=[\partial D]$.
\enddefinition

Shrink the reference fiber~$F$ to a closed disk containing
$\Conv_F\CC$ in its interior and $s(F)$ in its boundary. Pick a
geometric basis for~$\pi_F$ and identify it with a geometric basis
$\{\Ga_1,\Ga_2,\Ga_3\}$
for~$\AA$, establishing an isomorphism $\pi_F\cong\AA$.
Under this
isomorphism, the braid monodromy~$\bm$ takes values in the braid
group $\BG3\subset\Aut\AA$. The
monodromy~$\bm$ thus defined is independent of
the choice of a proper section, and another choice
of the geometric bases
for~$\pi_F$ and~$\AA$ results in the
global conjugation by a fixed braid $\Gb\in\BG3$, \ie, in the map
$\Gg\mapsto\Gb\1\bm(\Gg)\Gb$. Thus, the \emph{monodromy group}
$\BM_\CC:=\Im\bm\subset\BG3$
is determined by~$\CC$ up to conjugation.
One has $\depth\BM_\CC\divides|6d$;
the group $\BM_\CC$ is $\SS$-transitive if and only if $\CC$ is
irreducible.

Next statement is proved in~\cite{degt:trigonal}.

\theorem\label{th.monodromy}
The monodromy group of a non-isotrivial trigonal curve is of genus
zero. Conversely, given a subgroup $G\subset\BG3$ of genus zero
and depth $6d>0$, there is a unique, up to isomorphism and
$d$-Nagata equivalence, trigonal curve $\CC_G$ with the following
property\rom: for a non-isotrivial trigonal curve~$\CC$ one has
$\BM_\CC\prec G$ if and only if $\CC$ is $d$-Nagata equivalent to a
curve induced from~$\CC_G$.
This curve $\CC_G$
is called the
\emph{universal curve} corresponding to~$G$.
\qed
\endtheorem

The universal curve~$\CC_G$ can be reconstructed from the
skeleton~$\Sk_G$.
In fact,\mnote{more precise description}
$\Sk_G$ is the \emph{dessin d'enfants}, in the sense of Grothendieck, of a
unique (up to M\"{o}bius transformation of the source)
regular map $j\:\Cp1\to\Cp1=\C\cup\infty$ with three critical
values~$0$, $1$, and~$\infty$ only. This map~$j$ is the $j$-invariant
of~$\CC_G$ (thus defining~$\CC_G$ up to Nagata transformation),
and the types of the singular fibers
of~$\CC_G$ are given by the
type specification of~$G$ (which explains the term).

\subsection{Generalized curves}\label{s.generalized}
Now,\mnote{separate subsection; rewritten}
let $\CC\subset\Sigma_d$ be a generalized trigonal curve. This time, the
closure of $\Conv\CC$ does not need to be compact and $\CC$ may not admit a
proper section.
To overcome this difficulty, consider a proper model $\CC'\subset\Sigma_{d'}$
of~$\CC$ and, for a punctured disk $\disk^\circ$ as above,
denote by $\bm'\:\pi_1(\disk^\circ,F)\to\BG3$ the braid
monodromy of~$\CC'$.
Fix, further, a geometric basis
$\{\Gg_1,\ldots,\Gg_r\}$ for $\pi_1(\disk^\circ,F)$.
Then, the difference between~$\CC$ and~$\CC'$ can be described in terms of the
so-called \emph{slopes} $\slope_i\in\AA$ assigned to each geometric
generator~$\Gg_i$.
Roughly,\mnote{slopes defined}
assume that $\Gg_i$ is represented by a loop of the form
$l_i\cdot\mu_i\cdot l_i\1$, where $\mu_i$ is a small circle about a
fiber~$F_i$ and $l_i$ is a simple path connecting the
common base point and a point
$a_i\in\mu_i$.
Consider a small analytic disk $\Phi\subset\Sigma_d$ transversal to~$F_i$
and disjoint from~$\CC$ and~$E$, and a similar disk
$\Phi'\subset\Sigma_{d'}$ with respect to $\CC'$.
Let $\bar\Phi\subset\Sigma_{d'}$ be the transform of~$\Phi$, and assume
that the boundaries $\partial\Phi'$ and $\partial\bar\Phi$ have a common
point over~$a_i$.
Then, the loop $[\partial\bar\Phi]\cdot[\partial\Phi']\1$ is homotopic to a
certain class in the fiber over~$a_i$.
The image of this class
under the translation homomorphism along~$l_i\1$ is the slope;
it is
well defined up to a number of moves, irrelevant in the sequel.
For details and further
properties,
see \cite{degt:dessin}.

Now, the \emph{monodromy} of~$\CC$ is defined as the homomorphism
$\bm\:\Gg_i\mapsto\bm_i$, where
$\bm_i$ is the map $\Ga\mapsto\slope_i\1\bm_i'(\Ga)\slope_i$ and
$\bm_i'=\bm'(\Gg_i)$.
This monodromy takes values in the extended group $\Prod$;
its image $\BM_\CC$ is
called the \emph{monodromy group} of~$\CC$.
Strictly speaking,\mnote{the ambiguity discussed}
both~$\bm$ and $\BM_\CC$ depend on a number of choices
(trigonal model~$\CC'$, geometric basis~$\{\Gg_i\}$, slopes~$\slope_i$,
\etc.); however, we only retain the original curve~$\CC$ in the notation as
the other choices do not affect the fundamental group, \cf.
\autoref{th.vanKampen} below.

The projections $\pr_\MG\BM_\CC$ and $\pr_\MG\BM_{\CC'}$ coincide,
hence $\BM_\CC$ is also a subgroup of genus zero, see
\autoref{th.monodromy}. Unlike the case of genuine trigonal
curves, I do not know an intrinsic description of the subgroups of
$\Prod$ that can appear as the monodromy groups of generalized
trigonal curves.

\remark[Important remark]\label{BGvsBu}
It is worth emphasizing that the monodromy groups of genuine and
generalized trigonal curves lie, respectively, in the braid
group~$\BG3$ and extended group $\Prod$. Hence, all statements
below
concerning subgroups of $\Prod$ or $\Bu3$ hold for generalized
trigonal curves,
whereas those specific to subgroups of $\BG3$ hold for genuine
curves only. Formally, one can extend the statements concerning
subgroups $G\subset\BG3$ and \emph{extended} modules $\IM_G$ to
generalized trigonal curves with all slopes of degree divisible by
three.
\endremark

\subsection{The fundamental group}\label{s.pi1}
Consider a generalized trigonal curve $\CC\subset\Sigma_d$,
pick a
nonsingular fiber~$F_0$ of~$\CC$, and define the \emph{affine} and
\emph{projective
fundamental groups} of~$\CC$ to be
$\piaff\CC=\pi_1(\Sigma_d\sminus(\CC\cup E\cup F))$
and $\piproj\CC=\pi_1(\Sigma_d\sminus(\CC\cup E))$. The affine
group $\piaff\CC$ is an infinite cyclic central extension
of~$\piproj\CC$. In particular, the commutants of the two groups
are canonically isomorphic, hence so are the Alexander modules
defined below.

Fix all necessary data (trigonal model, proper section, bases,
an identification $\pi_F=\AA$, \etc., see
Subsections~\ref{s.sections} and~\ref{s.generalized})
and let $\BM_\CC$ be the resulting
monodromy group.
The following theorem is essentially contained
in~\cite{vanKampen}.

\theorem\label{th.vanKampen}
One has
$\piaff\CC=\AA/\<\Gb(\Ga)=\Ga,\ \Gb\in\BM_\CC,\ \Ga\in\AA\>$.
\qed
\endtheorem

It follows that $\piaff\CC$ depends on the conjugacy class of
$\BM_\CC\subset\Prod$ only.
Any presentation of $\piaff\CC$ as in \autoref{th.vanKampen}
is called \emph{geometric}. The group inherits from~$\AA$ the
degree homomorphism $\deg\:\piaff\CC\onto\Z$, which does not depend
on the choice of a geometric presentation. (The projective
group~$\piproj\CC$ is the quotient of~$\piaff\CC$ by a certain
central element of positive degree.)

Denote by~$\AM_\CC$ the abelianization of the kernel $\Ker\deg$. As
in Subsection~\ref{s.Burau}, the conjugation~$t$ by any element
$\Ga\in\piaff\CC$ of degree one turns~$\AM_\CC$ into a module
over~$\Lambda$; it is called the \emph{Alexander module} of~$\CC$,
and the order $\Delta_{\CC,p}\in\Lambda\otimes\Bbbk_p$ of the
$(\Lambda\otimes\Bbbk_p)$-module $\AM_\CC\otimes\Bbbk_p$, whenever
defined,
is called the
\emph{$({\bmod}\,p)$-Alexander polynomial} of~$\CC$.
In the classical setting, one usually considers
$\Delta_\CC:=\Delta_{\CC,0}$.
As an immediate consequence of \autoref{th.vanKampen}, one
concludes that $\AM_\CC=\AM_G$, where $G=\BM_\CC$ is the monodromy
group. For this reason, and in view of \autoref{th.monodromy},
in the rest of the paper we mainly deal with subgroups rather than
curves.

Letting $G=\BM_\CC$, one can also consider the extended module
$\bAM_\CC:=\bAM_G$, which `estimates' $\AM_\CC$ from above: there is
an epimorphism $\bAM_\CC\onto\AM_\CC$ (see \autoref{rem.eAP}).\mnote{ref
added}
Note however that
$\bAM_\CC$ is \emph{not} an invariant of the fundamental
group~$\piaff\CC$ only:
examples in Sections~\ref{S.proof.0} and~\ref{S.N=6}
show that $\bAM_\CC$ may be
nontrivial even when $\piaff\CC=\Z$.

\remark[Important remark]\label{important.remark}
Summarizing, one concludes that any upper bound on the
extended module~$\IM_G$ of a subgroup~$G$ of~$\Bu3$
(respectively,~$\BG3$)
of genus zero
can serve as an upper bound on the
conventional module $\tIM_\CC$ of a generalized (respectively,
genuine) trigonal curve~$\CC$. If $G$ is required to be
$\SS$-transitive, $\CC$ must be irreducible.
Furthermore,
according to \autoref{th.monodromy},
any finite index subgroup $G\subset\BG3$
of genus zero is the
monodromy group of a certain genuine trigonal curve. Hence, all
existence statements concerning subgroups of~$\BG3$ do imply the
existence of trigonal curves with desired properties.
\endremark

\remark\label{rem.meaning}
In view of \autoref{th.monodromy},
the
isomorphism $\AM_\CC=\AM_G$, $G=\BM_\CC$,
makes
Definitions~\ref{def.AP}--\ref{def.universal} geometrically
meaningful \emph{for subgroups of genus zero}. To generalize, one
could consider `trigonal curves' in geometrically ruled surfaces
$\Sigma\to\BB$ over arbitrary, not necessarily rational, bases.
However, in this case the presentation of $\piaff\CC$ is not
the one given by \autoref{th.vanKampen}: $\piaff\CC$ is the
quotient of the semidirect product
$\AA*\BM_C/\<\Gb\1\Ga\Gb=\Gb(\Ga),\ \Gb\in\BM_\CC,\ \Ga\in\AA\>$
by all elliptic and parabolic elements of $\BM_\CC$. (A subgroup
of genus zero is generated by its elliptic and parabolic elements,
see~\eqref{eq.G},
and one arrives at the statement of \autoref{th.vanKampen}.)
Thus, it is not quite
clear whether one should speak about the Alexander module of
$\piaff\CC$ itself (which is always large) or that of the kernel
of the inclusion epimorphism
$\piaff\CC\onto\pi_1(\Sigma)\cong\pi_1(\BB)$.
Nor is it clear how the universal subgroups should be defined in
this situation.
\endremark

\subsection{Plane curves with deep singularities}\label{s.plane}
Let $\DD\subset\Cp2$ be a plane curve with a distinguished
singular point~$P$ of multiplicity $\deg\DD-3$. Blow~$P$ up and
consider the proper transform~$\CC$ of~$\DD$: it is a generalized
trigonal curve in the Hirzebruch surface $\Sigma_1=\Cp2(P)$, the
exceptional section $E\subset\Sigma_1$ being the exceptional
divisor of the blow-up. The projection $\Sigma_1\to\Cp2$
establishes a diffeomorphism
\[*
\Sigma_1\sminus(\CC\cup E)\overset\cong\to\longrightarrow
 \Cp2\sminus\DD,
\]
hence an isomorphism $\piproj\CC=\pi_1(\Cp2\sminus\DD)$ of the
fundamental groups. Thus, all restrictions to the Alexander
module/polynomial of a generalized trigonal curve, in particular
Theorems~\ref{th.0}, \ref{th.p} and Addenda~\ref{add.N<=5},
\ref{add.N>10} in the introduction, hold
for plane
curves as above.
For this reason, we do not mention them separately.

\section{Local geometry of the skeleton\label{S.monodromy}}

In this section, we describe the local geometry of the skeleton of
a
finite index subgroup with nontrivial extended Alexander module.
The finite index condition is used in
Subsection~\ref{s.trivalent}: we assume that {\em all regions of
the skeleton are bounded}.

\subsection{Settings}\label{s.settings}
Fix a subgroup $G\subset\Bu3$ and let $\Sk=\Sk_G$ be its skeleton.
We assume that the index $[\MG:\bG]$ is finite, so that
$\Sk$ is a finite ribbon graph.

Fix, further, a field $\Bbbk=\Bbbk_p$ and an element
$\xi$ algebraic over~$\Bbbk$. Let $\bbbk=\Bbbk(\xi)$.
Unless stated otherwise,
we assume that $\xi\ne\pm1$.
Till the rest of the
paper, $M$ and~$N$ stand for the multiplicative orders of~$\xi$
and~$-\xi$, respectively. In particular we show that they are
finite.

In Subsections~\ref{s.trivalent}--\ref{s.white} below, we pick a
vertex~$v$
and an edge~$e$ close to~$v$, define a
certain subgroup $G_v\subset G$ generated by
some loops in a neighborhood of~$v$, and consider the
submodule
$\IM_v(\xi):=\IM_{G_v}(\xi)\subset\AM(\xi)$ and the
quotient $\bAM_v(\xi):=\AM(\xi)/\IM_v(\xi)$.
Then we introduce
a
basis $\{\Ga_1,\Ga_2,\Ga_3\}$ over~$e$,
see Subsection~\ref{ss.basis}, and use
this basis to analyze the conditions, `local' at~$v$,
necessary
for the nonvanishing $\bAM_v(\xi)\ne0$; the latter is
equivalent to the requirement that
$\dim_{\bbbk}\IM_v(\xi)\le1$ and is obviously necessary for
the nonvanishing $\bAM_G(\xi)\ne0$.

\subsection{A trivalent \pdfstr{black vertex}{\black-vertex}}\label{s.trivalent}
Consider a trivalent \black-vertex~$v$ of~$\Sk$; fix a marking~$e$
at~$v$ and a corresponding canonical basis
$\{\Ga_1,\Ga_2,\Ga_3\}$.
Let $G_v\subset G$ be the subgroup generated by the
boundaries of
$\reg(e)$ and $\reg(\X\1e)$,
\ie, by
$t^r\Gs_1^m$ and $t^s\Gs_2^n$, where
$m,n>0$ are the widths of the two regions and $r,s$ are given by
the corresponding type specifications.

Consider the matrix
$\CM=\bigl[\,t^{r}\Gs_1^m-\id\bigm|t^{s}\Gs_2^n-\id\,\bigr]$:
\[
\CM=\Mat{
 t^{r}(-t)^m-1&t^{r}\tf_m(-t)&t^{s}-1&0\\
 0&t^{r}-1&t^{s+1}\tf_n(-t)&t^{s}(-t)^n-1}.
\label{eq.M}
\]
Clearly, $\dim_{\bbbk}\IM_v(\xi)=\rank\CM(\xi)$, and we are
interested in the conditions on $m,n,r,s$ necessary and
sufficient for $\rank\CM(\xi)\le1$. Consider the following cases.

\subsubsection{Type~$\vZ$}\label{ss.minor.0}
If $\CM(\xi)=0$, \ie, $\IM_v(\xi)=0$,
the marked vertex~$v$ is said to be of \emph{type~$\vZ$}.
This is the case if and only if $N:=\ord(-\xi)<\infty$
divides both~$m$ and~$n$ and
$\ee_p(N)$ divides both~$r$ and~$s$.

Now, assume that $\IM_v(\xi)\ne0$ is a proper submodule of
$\AM(\xi)$.
Then $\xi$
annihilates all $(2\times2)$-minors of~$\CM$ and
one has one of the following three cases.

\subsubsection{Type~$\vI_1$}\label{ss.minor.1}
$\xi^{r}(-\xi)^m-1=\tf_m(-\xi)=\xi^{s}-1=0$, \ie,
the first row vanishes.
In this case, $N:=\ord(-\xi)\divides|m$; in addition, one has
$\xi^{r}=\xi^{s}=1$, \ie, $\ee_p(N)$ divides both~$r$
and~$s$.
The module~$\IM_v(\xi)$ is generated by $\tf_n(-t)\be_2$.
If $N\notdivides|n$,
the marked vertex~$v$ is
said to be of \emph{type~$\vI_1$}.
Then $\IM_v(\xi)\ne0$ is generated by~$\be_2$.

\subsubsection{Type~$\vI_2$}\label{ss.minor.2}
$\xi^{r}-1=\tf_n(-\xi)=\xi^{s}(-\xi)^n-1=0$, \ie,
the second row vanishes.
Similarly to the previous case,
$N:=\ord(-\xi)\divides|n$ and $\ee_p(N)$ divides both~$r$
and~$s$.
The module~$\IM_v(\xi)$ is generated by $\tf_m(-t)\be_1$.
If
$N\notdivides|m$,
the marked vertex~$v$ is
said to be of \emph{type~$\vI_2$}.
In this case, $\IM_v(\xi)\ne0$ is generated by~$\be_1$.

\subsubsection{Type~$\vII$}\label{ss.minor.2x2}
$\xi^{r}(-\xi)^m-1=\xi^{s}(-\xi)^n-1=\CM_{2,3}(\xi)=0$, where
$\CM_{2,3}$ is the minor composed of the second and third columns.
Modulo the first two relations,
$\xi^{r}\tf_m(-\xi)=(\xi^{r}-1)/(\xi+1)$ and
$\xi^{s+1}\tf_n(-\xi)=\xi(\xi^{s}-1)/(\xi+1)$,
see~\eqref{eq.identity},
and
$\CM_{2,3}(\xi)=-(\xi^{r}-1)(\xi^{s}-1)(\xi^2+\xi+1)/(\xi+1)^2$.
Thus, either
\roster
\item\label{2x2.1}
$\xi^{r}=1$, and then $N:=\ord(-\xi)\divides|m$, or
\item\label{2x2.2}
$\xi^{s}=1$, and then $N:=\ord(-\xi)\divides|n$, or
\item\label{2x2.0}
$\xi^2+\xi+1=0$.
\endroster

Using~\eqref{eq.identity} and the fact that $t+1$ is
invertible in~$\Lambda(\xi)$, one can see that the
module~$\IM_v(\xi)$ is generated by
$(t^{s}-1)((t\1+1)\be_1+\be_2)$ and
$(t^{r}-1)(\be_1+(t+1)\be_2)$
in Cases~\loccit{2x2.1} and~\loccit{2x2.2}, respectively.
In Case~\loccit{2x2.1}, assuming that $\xi^s\ne1$ (and hence
$N\notdivides|n$), the vertex is said to be of
\emph{type~$\vII_1$}; the module ~$\IM_v(\xi)$ is generated by
$(t\1+1)\be_1+\be_2$.
In Case~\loccit{2x2.2}, assuming that $\xi^r\ne1$ (and hence
$N\notdivides|m$), the vertex is said to be of
\emph{type~$\vII_2$}; the module ~$\IM_v(\xi)$ is generated by
$\be_1+(t+1)\be_2$.

In Case~\loccit{2x2.0}, assuming that $\xi^r\ne1$ and $\xi^s\ne1$,
we let $N=\ee_p(3)$ and assign to the vertex \emph{type~$\vIIex$}.
This is the only case when one cannot assert that $N\divides|m$ or
$N\divides|n$. (In fact, if $N$ does divide~$m$ or~$n$, then the
vertex is of type~$\vZ$, $\vII_1$, or~$\vII_2$.)
The module $\IM_v(\xi)$ is generated by any of the
two
elements $(t\1+1)\be_1+\be_2$ or $\be_1+(t+1)\be_2$ above.

Summarizing, one concludes that a necessary condition for the
nonvanishing
$\AM(\xi)/\IM_G(\xi)\ne0$ is that $N:=\ord(-\xi)<\infty$
and at each marked vertex $(v,e)$
other than of type~$\vIIex$ (which can only occur if $\ee_p(N)=3$)
at
least one of the regions
$\reg(e)$, $\reg(\X\1e)$
has width divisible by~$N$.

\definition\label{def.trivial}
With~$N$ fixed, a region of width divisible
by~$N$ is called \emph{trivial} (or \emph{$N$-trivial});
such a region does not contribute to $\IM_G(\xi)$.
A region of width not divisible by~$N$ is called
\emph{essential}, or \emph{$N$-essential}.
Essential regions are subdivided into type~$\vI$
and~$\vII$, depending on the type of the vertices in their
boundary.
\enddefinition

Summarizing, one arrives at the following statement.

\lemma\label{one.essential}
Assume that $\bAM_G(\xi)\ne0$ and let $M=\ee_p(N)$. Then\rom:
\roster
\item\label{1.trivial}
for each trivial region~$R$ one has
$\type(R)=\rdeg R\bmod2M$\rom;
\item\label{1.I}
for each type~$\vI$ essential region~$R$ one has
$\type(R)=\rdeg R\bmod2M$\rom;
\item\label{1.II}
for a type~$\vII$ essential region~$R$ of width $n=\rdeg R$,
if $n$ is even or $p=2$, then
$\type(R)=-n\bmod2M$, otherwise
$\type(R)=M-n\bmod2M$,
and in the latter case $M$ must be even\rom;
\item\label{1.<=1}
if $M\ne3$,
at each trivalent \black-vertex at most one region is
essential\rom;
\item\label{1.all}
if $\IM_G(\xi)=0$, then all regions are trivial.
\endroster
\par\removelastskip
\endlemma

\proof
Items~\loccit{1.trivial}--\loccit{1.II} paraphrase the
conditions $t^s=1$ and $t^s(-1)^n=1$ in terms of the type
specification.
For~\loccit{1.<=1}
and~\loccit{1.all}, it suffices to consider all three markings at the
given vertex or, respectively, at all vertices of the skeleton.
\endproof

\subsection{A monovalent \pdfstr{black vertex}{\black-vertex}
(type~\pdfstr{III}{$\vIII$})}\label{s.monovalent}
Consider a monovalent \black-vertex~$v$ and
let~$e$ be the adjacent edge.
In a canonical basis $\{\Ga_1,\Ga_2,\Ga_3\}$ over~$e$, the
positive loop about~$v$ lifts to an element of the form
$t^r(\Gs_2\Gs_1)$.
Let $G_v\subset G$ be the subgroup generated by this element.

One has $\det(t^r(\Gs_2\Gs_1)-\id)=\tf_3(t^{r+1})$.
Hence, one has
$\bAM_v(\xi)\ne0$ if and only if
$M:=\ord\xi<\infty$ satisfies the following conditions:
\roster*
\item
$M\divides|3(r+1)$ and $M\notdivides|(r+1)$
(in particular, $M=0\bmod3$)
if $p\ne3$, and
\item
$M\divides|(r+1)$ and $M\ne0\bmod3$ if $p=3$.
\endroster
If this is the case, the module~$\IM_v(\xi)$ is generated by
$-t^{r}\be_1+\be_2$.
Computing the exponents modulo~$M$, the latter can be rewritten in
the form $-t^{s}\be_1+\be_2$, where
\roster*
\item
$s=\pm\frac13M-1$ if $p\ne3$ and
\item
$s=-1$ (or $s=M-1$) if $p=3$.
\endroster
If $p\ne3$ then, according to the sign~$\pm$ in the expression
for~$s$ above, we assign to the vertex~$v$
\emph{type~$\vIII_\pm$}. If $p=3$, there is one type~$\vIII$.
Observe that, if $p\ne3$, the generator of $\IM_G(\xi)$ can be
rewritten in the form $-\epsilon t\1\be_1+\be_2$
with $\epsilon^2+\epsilon+1=0$.

Now, assume that $v$ has a trivalent neighbor~$u$ in~$\Sk$.
(The remaining cases are
treated in Subsection~\ref{s.special} below.)
Summarizing and using \autoref{one.essential}, one arrives at
the following statement.

\lemma\label{lem.black}
Assume that $\bAM_G(\xi)\ne0$ and let $M=\ee_p(N)$.
Let~$v$ be a monovalent \black-vertex and $u$ its trivalent
neighbor.
Then\rom:
\roster
\item\label{b.not3}
if $p\ne3$, then $M=0\bmod3$ and
$\type(v)=\pm\frac23M\bmod2M$\rom;
\item\label{b.p=3}
if $p=3$, then $M\ne0\bmod3$ and $\type(v)=0\bmod2M$\rom;
\item\label{b.trivial}
unless $u$ is of type~$\vIIex$, $v$ is in the boundary of an
$N$-trivial region.
\qed
\endroster
\par\removelastskip
\endlemma

\remark
If $G\in\BG3$ and $p\ne3$, see \autoref{lem.black}\loccit{b.not3},
the condition $\type(\BLACK)=2\bmod6$
in \autoref{type}\loccit{tp.black} implies that $M=\pm3\bmod9$
and, according to the sign in this congruence, only one
type~$\vIII_\pm$ can appear.
\endremark

\subsection{A monovalent \pdfstr{white vertex}{\white-vertex}
(type~\pdfstr{IV}{$\vIV$})}\label{s.white}
Consider a monovalent \white-vertex~$v$ and
let~$e'$ be the adjacent edge. To simplify the expressions below,
switch to the edge $e=\X\Y e'$.
In a canonical basis $\{\Ga_1,\Ga_2,\Ga_3\}$ over~$e$, the
positive loop about~$v$ lifts to an element of the form
$t^r(\Gs_2\Gs_1\Gs_2)$.
Let $G_v\subset G$ be the subgroup generated by this element.

One has $\det(t^r(\Gs_2\Gs_1\Gs_2)-\id)=1-t^{2r+3}$.
Hence, one has
$\bAM_v(\xi)\ne0$ if and only if
$M:=\ord\xi\divides|(2r+3)$,
and in this case $\IM_v(\xi)$ is generated by
$t^{r+1}\be_1+\be_2$,
which can be rewritten in the form
$t^s\be_1+\be_2$, where
$s=\frac12(M-1)$.

A monovalent \white-vertex~$v$ is said to be of
\emph{type~$\vIV$}. Assuming that $v$ is adjacent to a trivalent
vertex~$u$, one arrives at the following statement.

\lemma\label{lem.white}
Assume that $\bAM_G(\xi)\ne0$ and let $M=\ee_p(N)$.
Let~$v$ be a monovalent \white-vertex and $u$ its trivalent
neighbor.
Then\rom:
\roster
\item\label{w.M=1}
$M$ is odd and $\type(v)=M\bmod2M$\rom;
\item\label{w.trivial}
unless $u$ is of type~$\vIIex$, $v$ is in the boundary of an
$N$-trivial region.
\qed
\endroster
\par\removelastskip
\endlemma

\subsection{Two special subgroups}\label{s.special}
In this subsection, we treat the two cases
that are not quite covered by Lemmas~\ref{lem.black}
and~\ref{lem.white};
namely, we
consider a skeleton~$\Sk$ with a monovalent \black-- or
\white-vertex that is not adjacent to a trivalent \black-vertex.
Clearly, $\Sk$ is either
$\mathord\circ{\joinrel\relbar\joinrel\relbar\joinrel}\mathord\bullet$
or
$\mathord\bullet{\joinrel\relbar\joinrel\relbar\joinrel}\mathord\bullet$;
in the former case, $\bG=\MG$,
in the latter case, $\bG$
is the only
index~$2$ subgroup $\MG^2=2A^0$.

\proposition\label{th.MG}
If $\bG=\MG$, then $\AM_G(\xi)\ne0$ if and only if
$G\prec(\MG)\minus$.
In this case, one
has $p=2$, $\xi^2+\xi+1=0$, and
$\IM_G(\xi)=\bbbk(-t\be_1+\be_2)$.
\endproposition

\proof
It suffices to consider matrix~$\CM$ in~\eqref{eq.M} with $m=n=1$.
\endproof

\proposition\label{th.MG2}
Assume that $\bG=\MG^2=2A^0$ and $\AM_G(\xi)\ne0$.
Then either
\roster
\item\label{MG2.p=3}
$p=3$, $\xi=1$, and $\IM_G(\xi)=\bbbk(-t\be_1+\be_2)$\rom;
then $G\prec(2A^0)\two$, or
\item\label{MG2.p=0}
$\xi^2+\xi+1=0$ and
$\IM_G=\Lambda(-t\be_1+\be_2)\bmod\cp_3$\rom;
then $G\prec(2A^0)\minus$.
\endroster
\par\removelastskip
\endproposition

\proof
The group is generated by $t^r\Gs_1^2$, $t^s\Gs_2^2$,
$t^k\Gs_2\Gs_1$, and, possibly, an extra power of~$t$, and the
proof is a direct computation, starting
with~\eqref{eq.M} with $m=n=2$, \cf.
Subsections~\ref{s.trivalent} and~\ref{s.monovalent}.
\endproof

\remark\label{rem.MG}
The
largest subgroup of~$\MG$ (respectively, $\MG^2=2A^0$)
on which the slope $-{\bdeg}$ is equal to ${\bdeg}\bmod6$
is $\MG(3)=3D^0=\Ker({\bdeg}\bmod3)$
(respectively,
$\MG'=6A^1=\Ker({\bdeg}\bmod6)$;
this latter subgroup is of genus one).
\endremark

\subsection{A few consequences}\label{s.summary}
We state a few
immediate
consequences of the computation in
Subsections~\ref{s.trivalent}--\ref{s.white}. Note that in
\autoref{finite.order} we do \emph{not} assume that
$\depth G\ne0$ (which would make the claim trivial).

\lemma\label{finite.order}
If
$\bG\subset\MG$ is a subgroup of finite index,
there is an integer $M>0$ such
that $(t^M-1)(\AM_G\otimes\Bbbk_p)=0$ for each~$p$. In particular,
the Alexander polynomial $\bDelta_{G,p}$
is well defined and divides $(t^M-1)^2$.
\endlemma

\proof
One merely repeats the arguments of Subsections~\ref{s.trivalent}
and~\ref{s.special}, computing the ranks of the corresponding
matrices over $\Lambda\otimes\Bbbk_p$. Each time rank equals~$2$
and all invariant factors divide some $(t^M-1)$.
\endproof

\lemma\label{rank.2}
If $\IM_G(\xi)=0$, \ie, if $\rank\AM_G=2$,
then all vertices of~$\Sk_G$ are trivalent
\rom(equivalently, $G$ is torsion free\rom)
and all regions of~$\Sk_G$ are trivial.
\endlemma

\proof
According to Subsections~\ref{s.trivalent}, \ref{s.monovalent},
and~\ref{s.white}, each essential region of~$\Sk_G$ and each
monovalent vertex makes a nontrivial contribution to~$\IM_G(\xi)$.
\endproof

\subsection{Isotrivial curves}\label{s.isotrivial}
Recall that, in appropriate affine coordinates
$(x,y)$ in $\Sigma_d$, the
equation of an irreducible isotrivial genuine trigonal curve~$\CC$
can be written in the form
\[*
\textstyle
y^2=\prod_i(x-x_i)^{m_i},\qquad
 m:=\gcd(m_i)\ne0\bmod3.
\]
Hence, the monodromy group of any generalized trigonal curve
Nagata equivalent to~$C$ is the abelian group generated by
$t^r(\Gs_2\Gs_1)$ and $t^s\id$ for some $r,s\in\Z$.

\theorem
The extended Alexander polynomial of an irreducible isotrivial
generalized trigonal curve~$\CC$ divides $\tf_3(t^{r+1})$ for some
$r\in\Z$. If $\CC$ is a genuine curve, then
$\bDelta_{\CC,p}=\tf_3(t^{r+1})$ for some $r\in3\Z$ and any~$p$.
\endtheorem

\proof
Both statements follow from the description of the monodromy group
and the computation in Subsection~\ref{s.monovalent}. If $\CC$ is
a genuine curve, the monodromy group is generated by
$(\Gs_2\Gs_1)^m$, hence $r\in3\Z$ and $s=0$.
\endproof

\section{Proof of \autoref{th.p}\label{S.proof.p}}

Throughout this section, we fix~$p$ (a prime or zero), a subgroup
$G\subset\Bu3$ \emph{of genus zero}, and a root~$\xi$ of its
Alexander polynomial $\bDelta_{G,p}$. Let $N=\ord(-\xi)<\infty$.
Recall that we assume $\xi\ne\pm1$, hence $N\ge3$. The ultimate
goal of the section is a proof of \autoref{th.p} and the
estimate $N\le10$ for $p=0$, see \autoref{N<=10}.

\subsection{The boundary of a trivial region}\label{s.boundary}
Consider an $N$-trivial region~$R$ of a certain width~$Nm$.
With respect to the default marking, see
Subsection~\ref{ss.regions},
all vertices in~$\dR$ are of types~$\vZ$, $\vI_1$,
$\vII_1$, $\vIII_\pm$ (or~$\vIII$ if $p=3$), or~$\vIV$.
Define the \emph{distance} $\dist(v_1,v_2)\in\CG{Nm}$ between two
vertices $v_1,v_2\in\dR$ as the distance in~$R$, regarded as an
orbit of~$\X\Y$, between the corresponding edges $e_1$, $e_2$ used
in Subsections~\ref{s.trivalent}--\ref{s.white} to construct the
canonical bases.

\lemma\label{distance}
With two exceptions,
the distance in~$\dR$ between any two vertices of the same
type other than~$\vZ$ is divisible by~$N$. The exceptions are as
follows\rom:
\roster*
\item
$\ee_p(N)=3$ and the vertices are of type~$\vII$ or~$\vIII_-$, or
\item
$\ee_p(N)=3$, $p=2$, and the vertices are of type~$\vIV$.
\endroster
\par\removelastskip
\endlemma

\proof
Let $M=\ee_p(N)$.
Consider a vertex $v\in\dR$ of a type other than $\vZ$,
$\vI_2$,
or~$\vII_2$ (the two latter do not occur due to our choice of the
markings). According to
Subsections~\ref{s.trivalent}--\ref{s.white},
in
the corresponding canonical basis
the submodule~$\IM_v(\xi)$
is generated by a vector of the form $a_v(t)\be_1+\be_2$, where the
coefficient $a_v(t)\in\Lambda(\xi)$ depends on~$\xi$ and the type
of~$v$ only: one has
\[
a_v(t)=0,\quad t\1+1,\quad
 -t^s,\quad
 \text{or}\quad t^{\frac12(M-1)}
\label{eq.gens}
\]
for $v$ of type~$\vI_1$, $\vII_1$, $\vIII_\pm$ (or~$\vIII$ if $p=3$),
or~$\vIV$, respectively. Here, $s=\pm\frac13M-1$ for
type~$\vIII_\pm$ and $s=-1$ for type~$\vIII$.

Let $u\subset\dR$ be
another vertex at a distance~$d$ from~$v$.
Connecting the corresponding edges by a path in~$\dR$, one can
assume that the canonical bases used are related \via~$\Gs_1^d$,
and
a necessary condition for
$\bAM_G(\xi)\ne0$ is that the generators of~$\IM_u(\xi)$ and
$\Gs_1^d\IM_v(\xi)$ should be linearly dependent. This condition
results in the equation
\[
\tf_d(-\xi)\bigl((\xi+1)a_v(\xi)-1\bigr)=a_v(\xi)-a_u(\xi).
\label{eq.distance}
\]
If $u$ and~$v$ are of the same type, the right hand side vanishes
and~\eqref{eq.distance} takes the form $\tf_d(-\xi)=1$ or
$(\xi+1)a_v(\xi)=1$. In the former case, one has $N\divides|d$, as
stated; in the latter case, using the list above, one can see that
the equation either has no solutions (for type~$\vI_1$),
or implies $\xi=1$ (for type~$\vIII$ with $p=3$),
or implies $\xi^2+\xi+1=0$.
Indeed, if $v$ is a vertex of type~$\vII_1$, the equation
$(\xi+1)a_v(\xi)=1$ is equivalent to $\xi^2+\xi+1=0$. If $v$ is of
type~$\vIII_\pm$, then, switching to $a_v(t)=-\epsilon t\1$
with $\epsilon^2+\epsilon+1=0$, see Subsection~\ref{s.monovalent},
one has $\xi=-\epsilon/(\epsilon+1)$, hence $\xi^2+\xi+1=0$.
Furthermore, in this case $\epsilon\xi\1=\xi$, \ie, the type
is~$\vIII_-$.
If $v$ is of type~$\vIV$, then, letting
$s=\frac12(M-1)$ and hence $M=2s+1$, one has
\[*
[t^s(t+1)-1]\cdot t[t^s(t+1)+1]-[t^{2s+1}-1]\cdot(t+1)^2=t^2+t+1.
\]
Then $M=3$ and the equation turns into
$\xi^2+\xi=1$. Hence $p=2$.
\endproof

\lemma\label{I-II}
If $p=0$ and $N\ne4$, the boundary~$\dR$
cannot contain vertices of types both~$\vI_1$ and~$\vII_1$.
\endlemma

\proof
Assume that $v$ is of type~$\vI_1$ and $u$ is of type~$\vII_1$.
Then~\eqref{eq.distance} turns into the four term equation
$\xi(-\xi)^d+\xi^2+\xi+1=0$, in which each term is a root of
unity. Geometrically, the sum of four unit complex numbers equals
zero if and only if the summands split into two pairs of opposite
ones. Hence, the above equation implies $\xi=-1$ or $\xi=\pm i$;
in the latter case, one has $N=4$.
\endproof

\subsection{First estimates}\label{s.estimates}
Denote by~$R_i$ and~$S_j$, respectively, the trivial and essential
regions of~$\Sk$ (where $i$ and~$j$ run over certain index sets).
Introduce the following counts for~$\Sk$:
\roster*
\item
$v_1$ is the number of monovalent \black-vertices;
\item
$v_3$ is the number of trivalent \black-vertices;
\item
$e_1$ is the number of monovalent \white-vertices;
\item
$e_2$ is the number of edges connecting pairs of \black-vertices;
\item
$Nm_i$ is the width of the trivial region~$R_i$; let $m=\sum_im_i$;
\item
$n_j$ is the width of the essential region~$S_j$; let $n=\sum_jn_j$.
\endroster
For a trivial region~$R_i$, introduce also the following
parameters, counting special vertices in the boundary~$\dR_i$:
\roster*
\item
$\KI_i$ is the number of vertices of type~$\vI_1$ or~$\vII_1$;
\item
$\KIII_i$ is the number of vertices of type~$\vIII_\pm$
(or~$\vIII$ if $p=3$);
\item
$\KIV_i$ is the number of vertices of type~$\vIV$.
\endroster
For $*=\vI$, $\vIII$, or~$\vIV$, let $k^*_i=K^*_i/m_i$ and
$k^*=\max_ik^*_i$.
Unless $\ee_p(N)=3$,
in view of \autoref{distance}, one has
$0\le\kI,\kIII\le2$ and $0\le\kIV\le1$ and, due to \autoref{I-II}, one
has $\kI\le1$ if $p=0$ and $N\ne4$.
Furthermore, $\kIII\le1$ if $p=3$ and $\kIII$ and~$\kIV$ vanish
unless $\ee_p(N)$ satisfies certain divisibility conditions, see
Lemmas~\ref{lem.black} and~\ref{lem.white}
for the existence of vertices of the
corresponding types.

The total number of regions of~$\Sk$ does not exceed
$m+n$ and, since $\Sk$ is a ribbon graph of genus zero,
Euler's theorem implies
$m+n-e_2+v_1+v_3\ge2$. (The edges counted by~$e_1$ are cancelled
by  the monovalent \white-vertices.) As usual, one has
$v_1+3v_3=Nm+n=e_1+2e_2$, and, eliminating~$e_2$ and~$v_3$, one
can rewrite
Euler's inequality above in the form
\[*
(6-N)m+5n+4v_1+3e_1\ge12.
\]
Since all monovalent
vertices belong to the boundaries of trivial regions of~$\Sk$, see
Lemmas~\ref{lem.black} and~\ref{lem.white},
one has
\[*\textstyle
v_1=\sum_i\kIII_im_i\le\kIII m,\qquad
e_1=\sum_i\kIV_im_i\le\kIV m.
\]
Crucial is the following observation.

\lemma\label{d<m}
One has $n=\sum_i\kI_im_i\le\kI m$.
\endlemma

\proof
By definition, each vertex in the boundary of an essential region
is of type~$\vI$ or~$\vII$. On the other hand, due to
\autoref{one.essential}\loccit{1.<=1}, each such vertex~$v$ appears in the
boundary of an essential region exactly once (hence the number of
these vertices is~$n$) and admits a unique marking~$e$ with
respect to which it is of type~$\vI_1$ or~$\vII_1$. With this
marking,
$\reg(e)$
is a trivial region; hence $v$ is counted
exactly once in the sum $\sum_i\kI_im_i$.
\endproof

Substituting, one arrives at
\[\textstyle
(6-N)m+\sum_i(5\kI_i+4\kIII_i+3\kIV_i)m_i\ge12
\label{eq.N.exact}
\]
and
$(6+\kall-N)m\ge12$, where
$\kall=\max_i(5\kI_i+4\kIII_i+3\kIV_i)$.
This implies
\[
N<6+\kall\le6+5\kI+4\kIII+3\kIV.
\label{eq.N}
\]

\corollary\label{N<=26}
If $p=\fchar\Bbbk=0$, then $N\le21$\rom;
otherwise, $N\le26$.
\endcorollary

\proof
The statement follows from~\eqref{eq.N} and
the estimates on~$\kI$, $\kIII$, and~$\kIV$
listed right after their definition.
\endproof

\subsection{Further restrictions}\label{s.better.estimates}
We keep the notation introduced in Subsection~\ref{s.estimates}.

\table
\caption{Exceptional factors of~$\Delta$ (not realized)}\label{tab.others}
\bgroup
\def\*{\llap{$^{*}$}}%
\let\comma,\catcode`\,\active\def,{$\comma\ \ $}%
\def\get[#1]{$#1$\hss}%
\centerline{\vbox{\halign{\hss$#$&\quad\hss$#$\hss&\quad\get#\cr
\omit\hss$p$\hss&$N$&\omit\qquad
 Factors $\mp_\xi\in\FF{p}[t]$ of $\Delta$\hss\cr
\noalign{\vskip2pt \hrule}\cr
  3&13&[t^3+2t+1, t^3+2t^2+1, t^3+t^2+2t+1, t^3+2t^2+t+1]\cr
 23&11&[t+2, t+4, t+6, t+9, t+12, t+18]\cr
 29&14&[t+4, t+22]\cr
 31&15&[t+14, t+18, t+19, t+20]\cr
 37&12&[t+8, t+14, t+23, t+29]\cr
 43&14&[t+32, t+39]\cr
   &21&[t+14, t+40]\cr
 53&13&[t+28, t+36]\cr
 61&15&[t+16, t+42]\cr
 79&13&[t+38, t+52]\cr
127&21&[t+47, t+100]\cr
211&15&[t+83, t+150]\cr
\crcr}}}\egroup
\endtable

\lemma\label{sum.k}
If $N>10$ and the triple $(p,N,\mp_\xi)$ is not one of those
listed in Tables~\ref{tab.factors} and~\ref{tab.others},
then, for each trivial
region~$R_i$, one has $\kI_i+\kIII_i+\kIV_i\le1$.
\endlemma

\proof
It suffices to show that, under the assumptions, the distance
in~$\dR_i$ between any two vertices~$u$, $v$ of
types, respectively, $T_u$, $T_v$
other than~$\vZ$
is divisible by~$N$.
Due to \autoref{distance}, one can assume that
$T_u\ne T_v$.
Let $d=\dist(u,v)$. Then $\xi$ must
satisfy~\eqref{eq.distance} and, since the equation is obviously
$N$-periodic in~$d$, it suffices to consider the values
$d=1,\ldots,N-1$.

Now, for each $N=11,\ldots,26$, see \autoref{N<=26},
each $d=1,\ldots,N-1$, and each pair
$T_u\ne T_v$ of types,
consider the resultant~$\CR$ of~\eqref{eq.distance} and
$(-\xi)^N-1$. (All computations below were performed using
{\tt Maple}.) One has $\CR\ne0$, which proves the statement for
$p=0$. For each prime divisor $p\ne2,3$ of~$\CR$, consider the greatest
common divisor of the two polynomials over~$\FF{p}$, decompose it
into irreducible factors, and select those that do not divide
$(-\xi)^n-1$ for some $n<N$. The minimal polynomial $\mp_\xi$ must
be
one of these factors. The cases $p=2$ or~$3$ are
treated similarly, but
separately, as equation~\eqref{eq.distance} changes in these
cases.

The above procedure results in a finite collection (too large to
be listed here) of sequences
$(N,p,\mp_\xi;d,T_u,T_v)$.
For each triple $(N,p,\mp_\xi)$ thus
obtained, one can analyze the types of vertices that may appear
simultaneously in the boundary of a single region and improve the
\latin{a priory} estimate $\kall\le21$ used in~\eqref{eq.N}. (For
example, if all types that can appear in the same region
are~$\vI$, $\vII$, and~$\vIII_+$, the estimate improves to
$\kall\le14$, hence $N\le19$.)
Disregarding the triples that do not satisfy the
new inequality $N<6+\kall$, one obtains
Tables~\ref{tab.factors} and~\ref{tab.others}.
\endproof

\corollary\label{N<=10}
Unless $(p,N,\mp_\xi)$ is one of the triples
listed in Tables~\ref{tab.factors} and~\ref{tab.others},
one has $N\le10$.
\endcorollary

\proof
Replacing all coefficients in
the definition of~$\kall$ with their maximum~$5$
and using \autoref{sum.k}, one obtains
$\kall\le5$ in~\eqref{eq.N}.
\endproof

\subsection{Proof of \autoref{th.p} and \autoref{add.N>10}}\label{proof.p}
Due to \autoref{important.remark}, it suffices to
prove a similar statement for the extended Alexander modules
$\IM_G$ of subgroups $G\subset\Bu3$ of genus zero. Note that we do
not use $\SS$-transitivity.

Consider one of the triples $(p,N,\mp_\xi)$ listed in
Tables~\ref{tab.factors} and~\ref{tab.others}. The submodule
$\IM_G(\xi)$ has the form $\bbbk\bv$, where $\bv$ is one of
the vectors listed in~\eqref{eq.gens}.
We choose $\bv=\be_2$ and compute the genus of the corresponding
universal subgroup.
The
computation, using {\tt Maple}, proceeds as follows. Map $\BG3$
or~$\Bu3$ to
the finite group $\GL(2,\bbbk)$,
let $\CV=\bbbk\be_2$,
and enumerate the cosets modulo the universal
subgroup~$G_{\CV}$.
In order to pass to $\bG_{\CV}$, identify further $\CM$ and
$t^s\CM$ for $\CM\in\MAT(\bbbk)$ and $s\in\CG{M}$.
(If a subgroup of $\BG3$ is to be found, take only $s\in3\CG{M}$.)
The result
is the set of edges of the skeleton of~$G_{\CV}$, see
Subsection~\ref{s.Sk},
its \black-- and \white-vertices and regions being the orbits of
$\Gs_2\Gs_1$, $\Gs_2\Gs_1^2$,
and $\Gs_1$, respectively.
Compute the Euler characteristic and make sure that it equals~$2$.

After the computation is completed, one can use the cosets found
to verify that, in fact, all subspaces $\bbbk\bv$ with $\bv$
as in~\eqref{eq.gens} are conjugate to $\bbbk\be_2$; hence
they would produce the same universal subgroups.

This computation eliminates all triples listed in
\autoref{tab.others} and, for subgroups of~$\BG3$, the values
$p=5$ and~$19$ in \autoref{tab.factors}, thus completing the
proof of \autoref{th.p} and the existence part of
\autoref{add.N>10}. (For the existence, one should also use
\autoref{important.remark} and,
passing from~$\IM_G$ to~$\tIM_G$,
\autoref{cor.V=V}.)

Further analysis of the data obtained in \autoref{sum.k}
shows that, with~$N$
and~$p$ fixed, each triple $(d,T_u,T_v)$ gives rise to at most one
irreducible factor~$\mp_\xi$. Hence, this factor is uniquely
recovered from the geometry of any trivial region of the skeleton
containing vertices of more than one type (such a region must
exist to break the bound $\kall\le5$, \cf. \autoref{N<=10}),
and two distinct
factors cannot appear simultaneously.
\qed

\remark\label{coset.enum}
Since the images of~$\BG3$ and~$\Bu3$ in $\GL(2,\bbbk)$ are not
known, the coset enumeration procedure starts with
the identity and
keeps multiplying matrices by
$\Gs_2\Gs_1$ and $\Gs_1\Gs_2\Gs_1$,
comparing the result with all matrices already listed;
each
new matrix $\CM$
is added to the list together with all products
$t^s\CM$, $s=1,\ldots,M-1$.
(If $M=0\bmod3$ and a subgroup of~$\BG3$ is to be found, only
values $s=0\bmod3$ are used.)
The equivalence relation is
linear:
two matrices~$\CM_1,\CM_2\in\MAT(\bbbk)$
are equivalent if and only if $\bv^\perp(\CM_1-\CM_2)=0$.
\endremark

\example\label{example.7--10}
The elliptic case $N\le5$ and the parabolic case $N=6$ are treated
in details in Sections~\ref{S.proof.0} and~\ref{S.N=6} below,
while the range $7\le N\le10$ remains wide open. A few examples
are given in \autoref{tab.examples}; they were found by the
coset enumeration procedure described in Subsection~\ref{proof.p}
and \autoref{coset.enum}.
All groups listed are
$\SS$-transitive.
The notation is the same as in \autoref{tab.factors}: marked
with a~$^*$ are the triples that appear in the Alexander
polynomials of genuine trigonal curves, and the last column
describes the projection to~$\MG$ of the corresponding universal
subgroup. (This time, the universal subgroups in~$\BG3$ and
in~$\Bu3$ have the same projection to~$\MG$; they differ by the
depth.)

\table
\caption{Examples with $7\le N\le10$}\label{tab.examples}
\bgroup
\def\*{\llap{$^{*}$}}%
\let\comma,\catcode`\,\active\def,{$\comma\ \ $}%
\let\scolon;\catcode`\;\active\let;\comma
\def\get[#1]{$#1$\hss}%
\def\getG(#1;#2;#3;#4){$(#1\scolon#2\comma#3\scolon#4)$\hss}%
\centerline{\vbox{\halign{\hss$#$&\quad\hss$#$&\quad\get#&\quad\hss\getG#\hss\cr
\omit\hss$p$\hss&$N$&\omit\qquad
 Factors $\mp_\xi\in\Bbbk[t]$ of $\Delta$\hss&\omit\quad\hss$\bG\subset\MG$\hss\cr
\noalign{\vskip2pt \hrule}\cr
  2& \*7&[t^3+t+1, t^3+t^2+1]&(9;1;0;1^2 7^1)\cr
  3& \*8&[t^2+2t+2, t^2+t+2]&(10;0;1;1^2 8^1)\cr
  5& \*8&[t^2+2, t^2+3]&(78;0;0;1^6 8^9)\cr
 11&\*10&[t+2, t+6, t+7, t+8]&(24;2;0;1^2 2^1 10^2)\cr
 17& \*8&[t+2, t+8, t+9, t+15]&(36;0;0;1^4 8^4)\cr
 19&   9&[t+4, t+5, t+6, t+16, t+9, t+17]&(20;0;2;1^2 9^2)\cr
 29& \*7&[t+7, t+16, t+20, t+23, t+24, t+25]&(60;0;0;1^4 7^8)\cr
 37&   9&[t+7, t+9, t+12, t+16, t+33, t+34]&(76;0;4;1^4 9^8)\cr
 43& \*7&[t+4, t+11, t+16, t+21, t+35, t+41]&(132;0;0;1^6 7^{18})\cr
\crcr}}}\egroup
\endtable

Conjecturally,
Tables~\ref{tab.factors} and~\ref{tab.examples}
list \emph{all} triples $(p,N,\mp_\xi)$,
$N\ge7$ (including the case $p=0$),
that appear in the extended Alexander polynomials
of subgroups of~$\Bu3$, not necessarily
$\SS$-transitive, of genus zero.
The proof, in its current state, requires a
great deal of computation and a number of technical details still
need to be double checked. It will appear
elsewhere.
\endexample


\section{Proof of \autoref{th.0}\label{S.proof.0}}

In this section, we list all roots~$\xi$ of the extended Alexander
polynomials with $N:=\ord(-\xi)\le5$;
only finitely many universal subgroups appear, and they are all
congruence subgroups of genus zero.
Then we eliminate the remaining cases $6\le N\le10$ for $p=0$ and
prove \autoref{th.0}.

\subsection{Reduction to congruence subgroups}\label{s.cong}
For an integer $N\ge2$, denote by $\BG3(N)\subset\BG3$
the subgroup
normally generated by~$\Gs_1^N$. For $N\le5$ these subgroups are
of finite index: one has $\pr_\MG\BG3(N)=\MG(N)$ and
$\depth\BG3(N)=6$, $12$, $24$, and~$60$ for $N=2$, $3$, $4$,
and~$5$, respectively.

\lemma\label{action}
Fix an integer $N\ge2$ and denote $\AM'=\AM/\tf_N(-t)$. Then the
induced $\BG3$-action on $\AM'\rtimes\Z$ factors through
$\BG3/\BG3(N)$.
\endlemma

\lemma\label{G(N)}
In the notation of \autoref{action},
let $G_{\CV}$ be the universal subgroup corresponding to a
submodule $\CV\subset\AM'$. Then all cusp widths of~$G_{\CV}$
divide~$N$.
If $N\le5$, then $\bG_{\CV}\subset\MG$ is a congruence subgroup
of level $l\divides|N$.
\endlemma

\proof[Proof of Lemmas~\ref{action} and~\ref{G(N)}]
Any element conjugate to~$\Gs_1^N$ acts trivially
on~$\AM(\xi)$, see~\eqref{eq.powers}, and on the product
$\AM'\rtimes\Z$ (as $\Gs_1$ preserves~$\Ga_3$). If $N\le5$,
then $\MG(N)$ is normally generated by $\Bs_1^N$.
\endproof

\lemma\label{N<=5.faithful}
In the notation of \autoref{action},
if $N\le5$,
then
the action of the quotient $\BG3/\BG3(N)$ on
$\AM'\rtimes\Z$ is faithful.
\endlemma

\proof
For any $s\in\Z$ one has
\[
[(\Gs_2\Gs_1)^{3s}(\Ga_1)\cdot\Ga_1\1]=
 \tf_s(t^3)[(t-1)\be_1+(t^2-1)\be_2];
\label{eq.rho}
\]
hence the depth of the kernel of the action equals that
of~$\BG3(N)$, see above.

Specializing at $t=-1$, one obtains a faithful action of
$\SL(2,\CG{N})=\tMG/\tMG(N)$ on $\CG{N}\oplus\CG{N}$; hence, the
images of the kernel and of~$\BG3(N)$ in~$\MG$ also
coincide.
\endproof

The action of $\BG3/\BG3(6)$ is also faithful,
see \autoref{N=6.faithful},
but I do not know whether this statement extends to $N\ge7$.

Using Lemmas~\ref{action} and~\ref{G(N)} and the tables of
congruence subgroups found in~\cite{Cummins.Pauli}, one can easily
enumerate all conjugacy classes of submodules
$\IM_G\subset\AM/\tf_N(-t)$ for $N\le5$.
In Subsections~\ref{s.N=3}--\ref{s.many} below, we state a few
consequences in terms of the specializations
$\IM_G(\xi)\subset\AM(\xi)$.

\subsection{The cases $N=3$ and~$5$}\label{s.N=3}
In this subsection,
we do not assume \latin{a priori} that $G$ is of genus zero or
$\SS$-transitive. For each universal subgroup~$G_{\CV}$,
we indicate only
its image $\bG_{\CV}$; in each case,
the type specification is recovered uniquely
(sometimes up to automorphism) using \autoref{type}
(always $\depth G_{\CV}=2\ee_p(N)$)
and Lemmas~\ref{one.essential},
\ref{lem.black}, and~\ref{lem.white}.

\theorem\label{th.N=3}
Assume that the extended Alexander polynomial $\bDelta_{G,p}$ has
a root $\xi\in\bbbk\supset\Bbbk_p$, $p\ne3$, with $\ord(-\xi)=3$.
Then one has one of the following three mutually exclusive
cases\rom:
\roster
\item\label{N=3.dim=0}
$\IM_G(\xi)=0$\rom; then
$\bG\prec\MG(3)=3D^0$,
\autoref{fig.curves}\rom{(b)},
and unless $p=2$, one has
$\IM_G=0\bmod\cp_6$
\rom(hence $\IM_G(\xi)=0$ for any~$q$\rom)\rom;
\item\label{N=3.dim=1}
$\IM_G(\xi)\sim\bbbk\be_2$\rom; then
$\bG\prec\MG_1(3)=3B^0$,
\autoref{fig.curves}\rom{(a)},
and unless $p=2$, one has
$\IM_G\sim\Lambda\be_2\bmod\cp_6$
\rom(hence $\IM_G(\xi)\sim\bbbk_q\be_2$ for any~$q$\rom)\rom;
\item\label{N=3.p=2}
$p=2$ and $\IM_G(\xi)=\bbbk(-t\be_1+\be_2)$\rom; then
$G\prec(\MG)\minus$, see \autoref{th.MG}.
\endroster
If $G\subset\BG3$, then Case~\loccit{N=3.p=2} does not occur and
in all other cases
$2\tIM_G\subset\IM_G\bmod\cp_6$
and $\tIM_G(\xi)=\IM_G(\xi)$ unless $p=2$.
\endtheorem

The reason for the exception in Cases~\loccit{N=3.dim=0}
and~\loccit{N=3.dim=1} is the fact that, for $p\ne2$, the type
specification is defined modulo $2\ee_p(3)=12$, whereas for $p=2$ it
is only defined modulo $2\ee_2(3)=6$. Hence, the corresponding
universal groups are larger (index~$2$ extensions)
for $p=2$. The same remark applies to
\autoref{th.N=5} below.

\figure
\def\tl#1#2#3{\relax\hbox to0pt{\hss(#1)\enspace
$#3$\hss}}
\centerline{\vbox{\halign{\hss#\hss&&\qquad\hss#\hss\cr
 \cpic{z3}&\cpic{z3+z3}&\cpic{z5}&\cpic{z5+z5}\cr
 \noalign{\medskip}
 \tl a{}{3B^0\ (\cp_6)}&\tl b{}{3D^0\ (\cp_6^2)}&
 \tl c{}{5D^0\ (\cp_{10})}&\tl d{}{5H^0\ (\cp_{10}^2)}\cr}}}
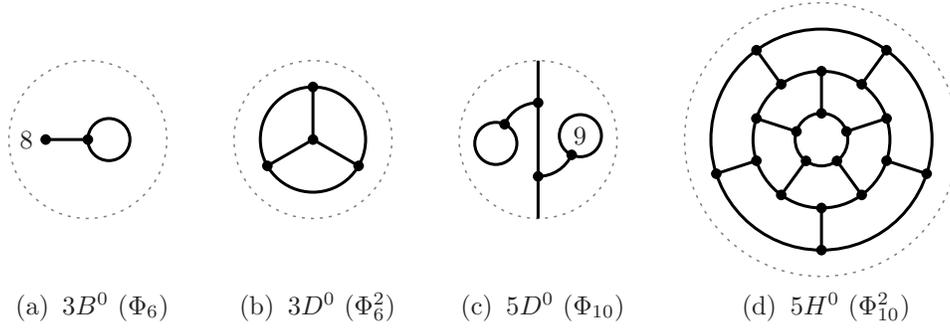
\caption{Skeletons of the universal groups for $p=0$}\label{fig.curves}
\endfigure

\proof
Assuming~$G$ universal and using \autoref{G(N)}, one concludes
that $\bG$ is a congruence
subgroup of level~$1$ or~$3$, and the submodules $\IM_G(\xi)$ can be
computed using the list found in~\cite{Cummins.Pauli}.
For $\bG=\MG(3)$ and $\MG_1(3)$,
one has $\IM_G=\cp_6\AM$ and
$\Lambda\be_2+\cp_6\AM$,
respectively, \cf.~\cite{degt:trigonal}.
The three other subgroups $3C^0\subset3A^0\subset\MG$
have $2$-torsion. Hence
$\IM_G(\xi)=\AM(\xi)$ unless $p=2$,
see \autoref{lem.white}\loccit{w.M=1}. If $p=2$, the universal
subgroup is given by \autoref{th.MG}.

Cases~\loccit{N=3.dim=1} and~\loccit{N=3.p=2} are mutually
exclusive since the largest subgroup of~$\MG_1(3)$ on which
the type specification shown in \autoref{fig.curves}(a) matches
$-{\bdeg}\bmod6$
is $\MG(3)$, see
\autoref{rem.MG}.
For the last statement, it suffices to notice that
the ideal generated by
$\cp_6$ and $t^2+t+1$
contains~$2\Lambda$,
hence $2\tIM_G\subset\IM_G\bmod\cp_6$,
see \autoref{cor.V}.
\endproof

\remark
If $p=2$ in \autoref{th.N=3}, the module $\tIM_G(\xi)$ depends
on the type specification, \ie, on the lift of~$G$, which can be
regarded as a subgroup of $\PSL(2,\FF3)=\MG/\MG(3)$,
to $\SL(2,\FF3)=\BG3/\BG3(3)$. In other words, $\IM_G(\xi)$
defines the type
specification
modulo~$6$, whereas
$\tIM_G(\xi)$ depends on its values modulo~$12$.
\endremark

\figure
\def\tl#1#2#3{\relax\hbox to0pt{\hss(#1)\enspace
$#3$\hss}}
\centerline{\vbox{\halign{\hss#\hss&&\qquad\hss#\hss\cr
 \cpic{4D0}&\cpic{5B0}&\cpic{5E0}&\cpic{5F0}\cr
 \noalign{\medskip}
 \tl a{}{4D^0}&\tl b{}{5B^0}&\tl c{}{5E^0}&\tl d{}{5F^0}\cr}}}
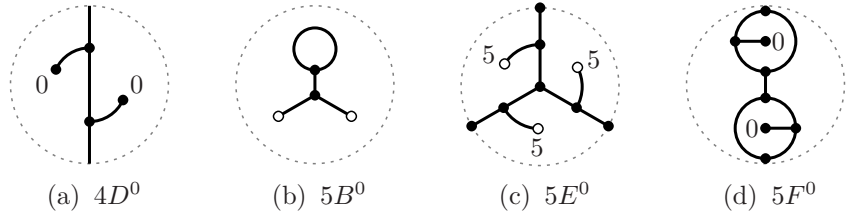
\caption{Some subgroups of level~$3$, $4$, and~$5$}\label{fig.groups}
\endfigure

\theorem\label{th.N=5}
Assume that the extended Alexander polynomial $\bDelta_{G,p}$ has
a root $\xi\in\bbbk\supset\Bbbk_p$, $p\ne5$, with $\ord(-\xi)=5$.
Then one has one of the following four mutually exclusive
cases\rom:
\roster
\item\label{N=5.dim=0}
$\IM_G(\xi)=0$\rom; then $\bG\prec\MG(5)=5H^0$,
\autoref{fig.curves}\rom{(d)},
and unless $p=2$, one has
$\IM_G=0\bmod\cp_{10}$
\rom(hence $\IM_G(\xi)=0$ for any~$q$\rom)\rom;
\item\label{N=5.dim=1}
$\IM_G(\xi)\sim\bbbk\be_2$\rom; then
$\bG\prec\MG_1(5)=5D^0$,
\autoref{fig.curves}\rom{(c)},
and unless $p=2$, one has
$\IM_G\sim\Lambda\be_2\bmod\cp_{10}$
\rom(hence $\IM_G(\xi)\sim\bbbk_q\be_2$ for any~$q$\rom)\rom;
\item\label{N=5.p=2}
$p=2$ and $\IM_G(\xi)\sim\bbbk(t^2\be_1+\be_2)$\rom;
then $\bG\prec5E^0$,
\autoref{fig.groups}\rom{(c)}\rom;
\item\label{N=5.p=3}
$p=3$ and $\IM_G(\xi)\sim\bbbk(\be_1-t\be_2)$\rom; then
$\bG\prec5F^0$, \autoref{fig.groups}\rom{(d)}.
\endroster
If $G\subset\BG3$, then $\tIM_G=\IM_G\bmod\cp_{10}$.
\endtheorem

\proof
As above, using \autoref{G(N)} one can assume that
$\bG\supset\MG(5)$ and use the list found in~\cite{Cummins.Pauli}.
The two torsion free subgroups $\bG=\MG(5)$ and $\MG_1(5)$
result in
$\IM_G=0\bmod\cp_{10}$ and $\Lambda\be_2\bmod\cp_{10}$,
respectively.
All other subgroups have torsion and, due to
Lemmas~\iref{lem.black}{b.not3} and~\iref{lem.white}{w.M=1},
one has
$\IM_G(\xi)=\AM(\xi)$ whenever $p\ne3$ and $G$ has $3$-torsion
or $p\ne2$ and $G$ has $2$-torsion.

Assume that $p=2$. The three
level~$5$ subgroups with $2$-torsion only
are $5B^0$ and $5G^0\subset5E^0$.
The skeleton of $5B^0$, see \autoref{fig.groups}(b),
contradicts \autoref{distance}. For the other two groups, a
direct computation shows that
$\IM_G(\xi)=\bbbk(t^2\be_1+\be_2)$.

Assume that $p=3$. The only level~$5$ subgroup with $3$-torsion
only is $5F^0$, see \autoref{fig.groups}(d). Over $\MG(5)$, it
is generated by $\Bs_2\Bs_1$ and,
lifting this element to $t^4(\Gs_2\Gs_1)$,
one obtains
$\IM_G(\xi)=\bbbk(\be_1-t\be_2)$, see
Subsection~\ref{s.monovalent}.

Cases~\loccit{N=5.dim=1}, \loccit{N=5.p=2}, and~\loccit{N=5.p=3}
are mutually
exclusive since the only common subconjugate of any pair of
corresponding universal subgroups is $\MG(5)$.

The last statement follows directly from \autoref{cor.V=V}.
\endproof

\subsection{The cases $N=1$, $2$, and~$4$}\label{s.N=1,2,4}
Here, we \emph{do} assume that the subgroup~$G$ is
$\SS$-transitive. Without this assumption, the number of cases in
Theorems~\ref{th.N=2}, \ref{th.N=1}, and~\ref{th.N=4} would be
much larger.
As above, we only indicate the image $\bG_{\CV}$ of the
universal subgroup~$G_{\CV}$;
the type specification is given by \autoref{type}
and Lemmas~\ref{one.essential},
\ref{lem.black}, and~\ref{lem.white}.

\theorem\label{th.N=2}
For an $\SS$-transitive subgroup~$G$,
assume that $\bDelta_{G,p}(1)=0$.
Then one has\rom:
\roster
\item\label{N=2.p=3}
$p=3$,
$\IM_G(1)=\bbbk(-t\be_1+\be_2)$, and
$G\prec(2A^0)\two$, see \autoref{th.MG2}\loccit{MG2.p=3}.
\endroster
The conventional Alexander polynomial $\Delta_{G,p}$ cannot vanish
at~$1$.
\endtheorem

\proof
Using \autoref{G(N)}, one can assume that
$\bG\supset\MG(2)$ and, for $G$ to be $\SS$-transitive,
$\bG$ must not lie in~$\MG_1(2)$, \ie, it must contain
$\MG^2=2A^0$, see~\cite{Cummins.Pauli}.
Thus, the statement about $\IM_G(1)$ follows from
Propositions~\ref{th.MG} and~\ref{th.MG2}, and a simple
computation for $\bG=2A^0$ (using the fact that
over $\MG(2)$ this subgroup is generated by
$\Bs_2\Bs_1$) shows that
$\tIM_G(1)=\AM(1)$.
\endproof

\theorem\label{th.N=1}
For an $\SS$-transitive subgroup~$G$,
assume that $\bDelta_{G,p}(-1)=0$.
Then one has one of the following five mutually exclusive cases\rom:
\roster
\item\label{N=1.p=3.0}
$p=3$ and $\IM_G(-1)=0$\rom; then
$\bG\prec\MG(3)=3D^0$,
\autoref{fig.curves}\rom{(b)}\rom;
\item\label{N=1.p=3.1}
$p=3$ and $\IM_G(-1)\sim\Bbbk\be_2$\rom; then
$\bG\prec\MG_1(3)=3B^0$,
\autoref{fig.curves}\rom{(a)}\rom;
\item\label{N=1.p=5.0}
$p=5$ and $\IM_G(-1)=0$\rom; then
$\bG\prec\MG(5)=5H^0$,
\autoref{fig.curves}\rom{(d)}\rom;
\item\label{N=1.p=5.1}
$p=5$ and $\IM_G(-1)\sim\Bbbk\be_2$\rom; then
$\bG\prec\MG_1(5)=5D^0$,
\autoref{fig.curves}\rom{(c)}\rom;
\item\label{N=1.p=7}
$p=7$ and $\IM_G(-1)\sim\Bbbk\be_2$\rom; then
$\bG\prec\MG_1(7)=7E^0$.
\endroster
If $G\subset\BG3$, then $\tIM_G=\IM_G\bmod(t+1)$.
\endtheorem

\proof
Essentially, the statement is the principal result
of~\cite{degt:trigonal}, where all
modules $\AM_G/(t+1)$
are classified. (Note that \autoref{G(N)} does not apply to
$N=1$.)
The action $\Bu3$ on $\AM/(t+1)$ factors
through~$\tMG$, and the
universal subgroups are of
the form $\tMG_m(n)$. There are five such subgroups of genus zero
that are $\SS$-transitive.
The last statement follows from \autoref{cor.V=V}.
\endproof

\theorem\label{th.N=4}
For an $\SS$-transitive subgroup $G\subset\Bu3$,
assume that $\bDelta_{G,p}$ has
a root $\xi\in\bbbk\supset\Bbbk_p$, $p\ne2$, with
$\ord(-\xi)=4$. Then one has\rom:
\roster
\item\label{N=4.p=3}
$p=3$,
$\IM_G(\xi)\sim\bbbk(\be_1-t\be_2)$, and
$\bG\prec4D^0$, \autoref{fig.groups}\rom{(a)}.
\endroster
If $G\subset\BG3$, then $\tIM_G=\IM_G\bmod(t^2+1)$.
\endtheorem

\proof
Using \autoref{G(N)}, one can assume that $\bG$ is a
congruence subgroup of level~$2$ or~$4$. According
to~\cite{Cummins.Pauli}, there are three $\SS$-transitive (\ie, not
subconjugate to $\MG_1(2)$) subgroups with this property:
$\MG^2=2A^0$, $4A^0$, and $4D^0$. All three have $3$-torsion;
hence $\IM_G(\xi)=\AM(\xi)$ unless $p=3$, see
\autoref{lem.black}\loccit{b.not3}. The subgroup~$2A^0$ was considered
in Subsection~\ref{s.special}. The subgroup~$4A^0$ has $2$-torsion
as well and is eliminated by \autoref{lem.white}\loccit{w.M=1}.
The
remaining subgroup $4D^0$, see \autoref{fig.groups}(a), is
generated over $\MG(4)$ by $\Bs_2\Bs_1$; lifting it to
$t^3(\Gs_2\Gs_1)$ and using Subsection~\ref{s.monovalent}, one
obtains $\IM_G(\xi)=\bbbk(\be_1-t\be_2)$.
The last statement follows
from \autoref{cor.V=V}.
\endproof

\subsection{Realizability and dependencies}\label{s.many}
We show that most pairs $(p,\mp_\xi)$ listed in the previous
two sections do appear in the
(extended)
Alexander polynomials of
genuine trigonal curves and that most of them are
mutually exclusive.

\theorem\label{th.BG}
With the exception of~\iref{th.N=3}{N=3.p=2}, each case listed in
Theorems \ref{th.N=3}, \ref{th.N=5}, \ref{th.N=2}, \ref{th.N=1},
and~\ref{th.N=4} can be realized by a subgroup of~$\BG3$,
\latin{ergo} by a genuine trigonal curve.
\endtheorem

\proof
The type specifications in \autoref{th.N=3}, except
Case~\iref{th.N=3}{N=3.p=2}, \emph{are} trivial modulo~$6$.
In all other theorems, one has $\gcd(M,3)=1$ and hence the
type specifications can be \emph{chosen} trivial modulo~$6$.
\endproof

\remark
The minimal, in the sense of the skeleton, genuine trigonal curve
with the Alexander polynomial~$\cp_{10}$ has non-simple
singularities. According to the type specification shown in
\autoref{fig.curves}(c), it must be a curve in~$\Sigma_{10}$
with the set of singular fibers
$\tJ_{8,0}\oplus2\tA_4\oplus\tA_0^*$.
\endremark

\theorem\label{th.cases}
The fourteen cases listed in
Theorems~\ref{th.N=3}, \ref{th.N=5}, \ref{th.N=2}, \ref{th.N=1},
and~\ref{th.N=4}, are related as follows\rom:
\roster
\item\label{<=>}
if $p\ne2$, then
\iref{th.N=3}{N=3.dim=0}$\;\Longrightarrow\;$\iref{th.N=1}{N=1.p=3.0}
and
\iref{th.N=3}{N=3.dim=1}$\;\Longrightarrow\;$\iref{th.N=1}{N=1.p=3.1}\rom;
\item\label{=>}
if $p\ne2$, then
\iref{th.N=5}{N=5.dim=0}$\;\Longrightarrow\;$\iref{th.N=1}{N=1.p=5.0}
and
\iref{th.N=5}{N=5.dim=1}$\;\Longrightarrow\;$\iref{th.N=1}{N=1.p=5.1}\rom;
\item\label{->}
\iref{th.N=4}{N=4.p=3}$\;\Longrightarrow\;$\iref{th.N=2}{N=2.p=3}\rom;
\item\label{<->}
Cases
\iref{th.N=3}{N=3.dim=1}$\;\Longrightarrow\;$\iref{th.N=1}{N=1.p=3.1}
and \iref{th.N=2}{N=2.p=3} can occur simultaneously\rom;
\item\label{<MG>}
Case~\iref{th.N=3}{N=3.p=2} can occur simultaneously with any
case except~\iref{th.N=3}{N=3.dim=0}, \ditto{N=3.dim=1}\rom;
\item\label{exclusive}
otherwise,
if $G$ is of genus zero,
the cases are mutually exclusive.
\endroster
The implications in~\loccit{<=>} turn into equivalences if
$G\subset\BG3$ and $p\ne2$.
\endtheorem

\proof
For~\loccit{<=>} and~\loccit{=>}, the universal subgroups
coincide, the type specifications in
Theorems~\ref{th.N=3} and~\ref{th.N=5}
being more restrictive (defined, respectively, modulo~$12$ or~$20$)
than those in
\autoref{th.N=1} (defined modulo~$4$ only).
If $G\subset\BG3$, \autoref{type} makes the
type specifications in \autoref{th.N=1}
well defined modulo $12$ as well (if $p\ne2$) and the
implications in~\loccit{<=>} turn into equivalences.

The implication in~\loccit{->} follows from the inclusion
$4D^0\subset2A^0$.

The fact that the cases within each theorem are mutually exclusive
is stated in the corresponding theorem.
Otherwise,
consider two cases and let $G_1$, $G_2$ be the corresponding
universal subgroups and $M_1$, $M_2$ the values of~$M$. For the
two cases to occur simultaneously, the projections $\bG_1$
and $\bG_2$ must have a common subconjugate of genus zero.
Then, if in addition $\gcd(M_1,M_2)=1$, one can
also find a common type
specification.
Common subconjugates can be analyzed using
the tables found in~\cite{Cummins.Pauli} (listing, in
particular, all sub-/supergroups). Unless one of the groups
is~$\MG$ itself (Item~\loccit{<MG>} of the statement), the only
pair is
$2A^0\supset6C^0\subset3B^0$,
which accounts for Item~\loccit{<->}.
\endproof

\corollary
For an irreducible genuine trigonal curve~$C$, if
$\cp_M^2\divides|\Delta_C$, then the fundamental group $\piaff{C}$
admits a dihedral quotient $\GDG(\CG{M}\oplus\CG{M})$.
\qed
\endcorollary

This fact was stated in~\cite{degt:trigonal} without proof.

\theorem\label{th.rk=2}
For an $\SS$-transitive subgroup $G\subset\Bu3$ of genus zero,
assume that $\IM_G(\xi)=0$ for some
$\xi\in\bbbk\supset\Bbbk_p$. Then
$\bG\subset\MG(N)$ for $N=3$ or~$5$, see
Cases \iref{th.N=3}{N=3.dim=0}, \iref{th.N=5}{N=5.dim=0},
and~\iref{th.N=1}{N=1.p=3.0} and~\ditto{N=1.p=5.0}.
\endtheorem

\proof
According to \autoref{rank.2}, $G$ is torsion free and all its
cusp widths are divisible by~$N$. As a torsion free subgroup of
genus zero, $G$ is generated by its parabolic elements, hence
$G\subset\MG(N)$. Then $N\le5$, and it remains to observe that the
subgroups $\MG(4)\subset\MG(2)$ are not $\SS$-transitive.
\endproof

\subsection{The case $p=0$ and
\pdfstr{e\sb0(N)=q\textcircumflex r}{$\ee_0(N)=q^r$}}\label{s.M=q}
Assume that one of the polynomials $\Delta_G$ or
$\bDelta_G$ of a
subgroup $G\subset\Prod$ has a root~$\xi$ of order $M:=\ee_0(N)=q^r$,
where $q$ is a prime. In the range $3\le N\le10$, see
\autoref{N<=10}, this is the case for $N=4$, $6$, $8$,
or~$10$, \ie, for all even values of~$N$.

\lemma\label{prime.power}
Let $M=q^r$ be a prime power.
If $\cp_M$ divides $\Delta_G$ or $\bDelta_G$, then
$(t-1)$ divides $\Delta_{G,q}$ or~$\bDelta_{G,q}$, respectively.
\endlemma

\proof
We will prove the statement for~$\Delta$; the proof
for~$\bDelta$ is a literal repetition.

Under the assumptions, the group $\AM_G/\cp_M$ is
infinite (as it remains nontrivial after tensoring with~$\Q$).
Hence, the $q$-group $\Hom(\AM_G/\cp_M,\FF{q})$ is nontrivial, and
the order $q^r$ automorphism~$t$ of this group has a nontrivial
invariant element~$\Gf$. Then the $(\Lambda\otimes\FF{q})$-module
$\AM_G\otimes\FF{q}$ factors to
$\Im\Gf\cong\FF{q}=(\Lambda\otimes\FF{q})/(t-1)$.
\endproof

\corollary\label{th.M=q}
Let $M=q^r$ be a prime power, and let $G$ be an
$\SS$-transitive subgroup.
Then $\cp_M\notdivides|\Delta_G$, and if
$\cp_M\divides|\bDelta_G$, one has\rom:
\roster
\item
$M=3$, $\IM_G=\Lambda(-t\be_1+\be_2)\bmod\cp_3$ and
$G\prec(2A^0)\minus$, see
\autoref{th.MG2}.
\endroster
A subgroup~$G$ with these properties cannot lie in~$\BG3$.
\endcorollary

\proof
The statement follows from \autoref{prime.power} and
Theorems~\ref{th.N=2},~\ref{th.N=6}.
\endproof

\subsection{Eliminating $N=7$ and~$9$ for $p=0$}\label{s.N=7,9}
Let $\Bbbk=\Q$ and let $\xi$ be a primitive root of~$(-1)$ of
degree~$7$ or~$9$, so that $\mp_\xi$ is the cyclotomic polynomial
$\cp_{14}$ or $\cp_{18}$,
respectively. Note that, in both cases, $\deg\mp_\xi=6$, so that
$\Q(\xi)\supset\Q$ is a Galois extension of degree six.

Fix a vector $h\in\AM(\xi)$ and consider the universal subgroup
\[*
G_h:=\{\Gb\in\Bu3\,|\,\Im[\Gb(\xi)-\id]\subset\Lambda(\xi)h\}
 \subset\Bu3.
\]

\lemma\label{lem.Burau.infinite}
For any $h\in\AM(\xi)$, one has $[\Bu3:G_h]=[\MG:\bG_h]=\infty$.
\endlemma

\proof
Consider the element
\[*
\Gb:=t\Gs_1\1\Gs_2=\mat[t-1, -t, t^2, -t^2]
\]
and its specialization $\Gb(\xi)$.
We assert that, in an appropriate extension
$\bbbk\supset\Q(\xi)$ of degree at most two,
$\Gb(\xi)$
has two distinct eigenvalues {\em which are not roots of unity.}
Indeed, the characteristic polynomial of $\Gb(\xi)$ is
$\chi(\Gl)=\Gl^2+(\xi^2-\xi+1)\Gl+\xi^2$ and its
roots belong to an extension of~$\Q$ of
degree~$6$ or~$12$.
Hence, the degree of the minimal polynomial $\mp_\Gl\in\Q[t]$ of
any eigenvalue~$\Gl$ divides~$12$.
There are finitely many cyclotomic polynomials
$\cp_n$ with $\deg\cp_n\divides|12$ (one has
$n=1$, $2$, $3$,
$4$,
$5$, $6$, $7$, $8$, $9$, $10$, $12$, $13$, $14$, $18$,
$21$,
$26$, $28$, $36$,
or~$42$;
alternatively,
in the computation below
one can use the
polynomials $\Gl^n-1$ with $n=8$, $10$, $26$, $28$, $36$,
or~$42$). For
each such polynomial $\cp_n(\Gl)$, compute the resultant
$\CR_n(\xi)$ of $\cp_n(\Gl)$ and $\chi(\Gl)$ with respect
to~$\Gl$, treating~$\xi$ as an independent variable. Each time, it
is straightforward that $\CR_n\bmod\mp_\xi\ne0$; hence $\chi(\Gl)$
and $\cp_n(\Gl)$ have no common roots in any extension of
$\Q(\xi)$. (This computation was performed using {\tt Maple}.) It
follows that the two roots of~$\chi$ are not roots of unity and, in
particular, they are distinct (as their product~$\xi^2$ is a root
of unity).

Thus, for any pair $m\ne0$, $r$ of integers, the two eigenvalues
of $\xi^r\Gb^m(\xi)-\id$ are distinct and both nonzero. Hence,
$\rank[\xi^r\Gb^m(\xi)-\id]=2$ and $t^r\Gb^m\notin G_h$. On the
other hand, the projection $\Bb\in\MG$ is an element of infinite
order.
\endproof

\corollary\label{no7,9}
If $G\subset\Bu3$
and $[\MG:\bG]<\infty$ \rom(\eg, if $G$ is a
subgroup of genus zero\rom),
the
polynomial $\bDelta_G$ is not divisible by $\cp_{14}$ or
$\cp_{18}$.
\qed
\endcorollary

\remark
In the proof\mnote{new remark}
of \autoref{lem.Burau.infinite}, we used \Maple\ to show that a certain
algebraic number is not a root of unity. Probably, there should be a better
way to detect rational arguments, and I expect that the statement of the lemma
holds for any primitive root $\xi\in\C$ of $(-1)$ of degree $N\ge7$.
\endremark

\subsection{Proof of \autoref{th.0}}\label{proof.0}
Let $G\subset\Prod$ be the monodromy group of the curve; it is an
$\SS$-transitive subgroup of genus zero, see
\autoref{th.monodromy}.
Due to \autoref{finite.order} and \autoref{N<=10}, each
irreducible factor of~$\Delta_G$ is of the form~$\cp_{M}$ with
$N:=\ee_0(M)\le10$. Most values of~$N$ are eliminated above, see
\autoref{th.N=2} for $N=2$, \autoref{th.N=1} for $N=1$,
\autoref{th.M=q} for $N=4$, $6$, $8$, and~$10$, and
\autoref{no7,9} for $N=7,9$.
The multiplicity of each of the remaining factors~$\cp_6$,
$\cp_{10}$ cannot exceed two,
see \autoref{finite.order},
and two distinct factors cannot appear
simultaneously according to \autoref{th.cases}.
The realizability is given by \autoref{th.BG},
see \autoref{important.remark}.
\qed

\section{The case $N=6$\label{S.N=6}}

In this section, we treat the parabolic case $N:=\ord(-\xi)=6$.
Since $\xi^2+\xi+1=0$ in this case, \autoref{cor.V=V} does
not apply and we consider the extended modules
$\IM_G\subset\AM$ and $\bAM_G$ only.

\subsection{The action on~$\AM'$}\label{s.action}
Let $\Lambda'=\Lambda/\cp_3$ and $\AM'=\AM/\cp_3$; for an integer
$m>1$, let also $\Lambda'_m=\Lambda'\otimes\CG{m}$ and
$\AM'_m=\AM'\otimes\CG{m}$. Consider the vector
$\bv=-t\be_1+\be_2\in\AM'$. It is immediate that
$\Gs_1(\bv)=\Gs_2(\bv)=\bv$, and in the basis $\{\bv,\be_2\}$ the
induced $\BG3$-action is given by the matrices
\[
\Gs_1=\mat[1,-t^2,0,-t],\quad
\Gs_2=\mat[1,0,0,-t];
\label{eq.N=6}
\]
hence
\[
\Gs_1\Gs_2\1=\mat[1,t,0,1],\quad
\Gs_2\1\Gs_1=\mat[1,t+1,0,1].
\label{eq.MG'}
\]
It follows that the image of the action
on~$\AM'$ is the full group of upper
triangular matrices with $[1,(-t)^s]$ in the diagonal.

Let $\MG'=[\MG,\MG]=6A^1$; recall that it is the free
subgroup generated by $\Bs_1\Bs_2\1$ and
$\Bs_2\1\Bs_1$.
Let further $\MG''=[\MG',\MG']$ be the second commutant,
and let
$\MG''_m\subset\MG'$ be the preimage of $m\Z\oplus m\Z$ under the
abelianization homomorphism $\MG'\to\Z\oplus\Z$.

Next two lemmas follow immediately from~\eqref{eq.N=6}
and~\eqref{eq.MG'}.

\lemma\label{N=6.kernel}
The kernels of the $\BG3$-actions on~$\AM'$ and~$\AM'_m$ are the
subgroups $(\MG'')\plus$ and $(\MG''_m)\plus$, respectively. The
image of $(\MG')\plus\!/(\MG'')\plus$ in $\SL(\AM')$ consists of
all unipotent upper triangular matrices.
\qed
\endlemma

\lemma\label{conjugate}
Any two vectors of the form $\be_2+f_i\bv\in\AM'$,
$f_i\in\Lambda'$, $i=1,2$,
are conjugate to each other.
\qed
\endlemma

\lemma\label{N=6.faithful}
The action of $\BG3/\BG3(6)$ on~$\AM'\rtimes\Z$, see
\autoref{action}, is faithful.
\endlemma

\proof
The images in~$\MG$ of the elements $\Gb_1:=\Gs_1\Gs_2\1$ and
$\Gb_2:=\Gs_2\1\Gs_1$ generate~$\MG'$, and the image of the commutator
$[\Gb_1,\Gb_2]:=\Gb_1\Gb_2\Gb_1\1\Gb_2\1$ normally generates
$\MG''$.
One can easily check the identity
$[\Gb_1,\Gb_2](\Gs_2\Gs_1)^{-3}=\Gs_1\Gs_2^{-6}\Gs_1\1$; hence
$\BG3(6)$ is the lift $(\MG'')^0$ of $\MG''$ with the slope
$0\:\MG''\to\Z$.
On the other hand,
due to~\eqref{eq.rho}, the kernel of the action is a subgroup
of depth~$0$.
\endproof

\lemma\label{transitive}
If $m\ne0\bmod3$, the subgroup $\MG''_m$ is $\SS$-transitive.
\endlemma

\proof
One has $(\Bs_1\Bs_2\1)^m\in\MG''_m$.
\endproof

\lemma\label{genus}
Any subgroup $G\subset\MG$ containing~$\MG''_m$
is of genus at most one. If $m$ is
prime to~$6$, the following statements are equivalent\rom:
\roster
\item\label{g.zero}
$G$ is of genus zero\rom;
\item\label{g.index}
$[\MG:G]\ne0\bmod6$\rom;
\item\label{g.MG}
$G\not\subset\MG'$\rom;
\item\label{g.torsion}
$G$ has torsion.
\endroster
\par\removelastskip
\endlemma

\proof
The group~$\MG'$ is torsion free and all its cusp widths are equal
to~$6$. According to \autoref{G(N)}, all cusp widths
of~$\MG''_m$ divide~$6$. Hence, the covering
$\Sk_{\MG''_m}\to\Sk_{\MG'}$ is unramified,
see Subsection~\ref{ss.coverings},
and $\MG''_m$ is of
genus one (as $\MG'$ is of genus one and
any unramified covering of a torus is a torus).

One has $[\MG:\MG']=6$ and $[\MG':G']\divides|m^2$ for any
$G'\subset\MG'$; hence, statements~\loccit{g.index}
and~\loccit{g.MG} are equivalent. Obviously, \loccit{g.index}
implies~\loccit{g.torsion}, and \loccit{g.torsion}
implies~\loccit{g.MG}, as $\MG'$ is torsion free.
Since $\MG'$ is of genus one, \loccit{g.zero}
implies~\loccit{g.MG}. Finally, if $G$ has torsion, the covering
$\Sk_{G\cap\MG'}\to\Sk_G$ is ramified, see
Subsection~\ref{ss.coverings},
and $G$ is of genus zero (as
any \emph{ramified} covering by a torus has sphere as the base);
thus, \loccit{g.torsion} implies~\loccit{g.zero}.
\endproof

\subsection{Subgroups of $\BG3$}\label{s.trivial}
In this and next subsections, we treat the case of genuine
trigonal curves, \ie, we assume that $G\subset\BG3$.

\lemma\label{module}
Assume that $m$ is prime to~$6$. For a subgroup
$G\subset\BG3$, denote by $\IM'_G$ the projection of~$\IM_G$
to~$\AM'_m$. Then either
\roster
\item\label{m.0}
$\IM'_G=0\bmod\Lambda'_m\bv$, and then
$\bG\subset\MG'$, or
\item\label{m.all}
$\IM'_G=\AM'_m\bmod\Lambda'_m\bv$, and then
$\bG\not\subset\MG'$.
\endroster
Conversely, any submodule $\CV\subset\AM'_m$
satisfying~\loccit{m.0} or~\loccit{m.all} above is of the form
$\IM'_G$ for some subgroup $G\subset\BG3$.
\endlemma

\proof
All statements follow
immediately from the description of the
action \via\ upper triangular matrices, see~\eqref{eq.N=6}
and~\eqref{eq.MG'},
and the fact that
all polynomials $(-t)^s-1$, $s\ne0\bmod6$,
are invertible in $\Lambda'_m$.
\endproof

\theorem\label{N=6.BG}
Let $G\subset\BG3$ be a subgroup of genus zero. Then the module
$\AM'_G:=\AM_G/\cp_3$ is finite and, modulo $2$- and $3$-torsion,
one has
\roster
\item
$\AM'_G=\AM'_m/(\Lambda'_m\be_2+I\bv)$
for some integer $m$ prime to~$6$ and
ideal $I\subset\Lambda'_m$.
\endroster
Conversely,
any module
$\AM'_m/(\Lambda'_m\be_2+I\bv)$ as above is of the form
$\AM'_G$ for some $\SS$-transitive subgroup $G\subset\BG3$ of
genus zero.
\endtheorem

\proof
One has $(\Gs_1\Gs_2\1)^s\be_2-\be_2=st\bv$ for any $s\in\Z$.
Hence, for $\AM'_G$ to be infinite, the submodule
$\IM_G/\cp_3$ must lie in $\Lambda'\bv$. Then
\autoref{module} implies that $\bG\subset\MG'$ is
a subgroup of genus at lest one.

Assume that $\AM'_G$ is finite. Then, modulo $2$- and $3$-torsion,
$\AM'_G=\AM'_m/\IM'_G$ for some sufficiently large $m$ prime
to~$6$. Since $\bG\not\subset\MG'$, one has
$\IM'_G=\AM'_m\bmod\Lambda'_m\bv$, see \autoref{module}, \ie,
$\IM'_G$ contains a vector of the form $\be_2+f\bv$,
$f\in\Lambda'_m$. In view of \autoref{conjugate}, any such vector
is conjugate to~$\be_2$, \ie, up to conjugation
$\IM'_G$ is as stated in the theorem.

Conversely, any submodule $\CV=\Lambda'_m\be_2+I\bv\subset\AM'_m$ as
in the statement is of the form $\IM'_G$ for some subgroup
$G\subset\BG3$ with $\MG''_m\subset\bG\not\subset\MG'$, see
\autoref{module}; this subgroup is $\SS$-transitive,
\autoref{transitive}, and of genus zero, \autoref{genus}.
\endproof

\subsection{A characterization of universal subgroups}\label{s.universal}
Let $\Sk$ be the skeleton of a genus zero subgroup $G\subset\MG$.
Assume that it has $1$-, $2$-, $3$-, and
$6$-gonal regions only. Then Euler's formula yields
\[
3\nwhite+4\nblack+5n_1+4n_2+3n_3=12,
\label{eq.Euler}
\]
where $\nwhite$, $\nblack$, and~$n_i$, $i=1,2,3,6$ are the numbers
of, respectively, monovalent \white-- and \black-vertices and
$i$-gonal regions of~$\Sk$. One of the solutions to this equation
is $\nwhite=\nblack=n_1=1$, $n_2=n_3=0$, and in this case one has
$[\MG:G]=6n_6+1$.

\definition
A proper finite index subgroup $G$ of~$\BG3$ (or of~$\MG$)
is called \emph{$6$-significant} if
$\depth G=6$ and the skeleton~$\Sk$ of~$G$
has exactly two monovalent vertices,
one \white-- and one \black--, and one monogonal region, while all
other regions of~$\Sk$ are hexagons.
\enddefinition

Note that any $6$-significant subgroup $G\subset\BG3$ is automatically
$\SS$-transitive and of genus
zero.
Note also that, since $\depth G=6$, it is not important
whether we speak about subgroups of~$\BG3$ or~$\MG$:
one always has $G=(\bG)\plus$.
Examples of $6$-significant subgroups are shown in
\autoref{fig.N=6}.

\lemma\label{significant}
Any $6$-significant subgroup $G\subset\MG$ contains $\MG''_m$ for
some integer $m$ prime to~$6$. One can take $m=[\MG:G]$.
\endlemma

\proof
Since $[\MG:G]=6n_6+1=1\bmod6$, see above, the subgroup
$G':=G\cap\MG'$ is of index~$6$ in~$G$, torsion free, and with all
cusp widths $0\bmod6$. On the other hand, there is a unique
subgroup of~$G$ with these properties: it corresponds to the
$6$-fold cyclic covering
${\Sk'}\to\Sk_G$ appropriately ramified over the
monovalent vertices and the monogonal region. From the latter
description, it follows that all cusp widths of~$G'$ are
\emph{equal} to~$6$. Hence, $\Sk_{G'}\to\Sk_{\MG'}$ is an
unramified covering of degree
$m=[\MG':G']=[\MG:G]$ and, the fundamental
group of the torus being
$\Z\times\Z$,
it splits the
$(\CG{m}\times\CG{m})$-covering corresponding to
the inclusion $\MG''_m\into\MG'$.
\endproof

\theorem\label{th.significant}
A subgroup $G\subset\BG3$ is the universal subgroup corresponding
to a proper submodule $\CV=\Lambda'_m\be_2+I\bv\subset\AM'_m$ for some
sufficiently large~$m$ prime to~$6$, see \autoref{N=6.BG},
if and only if it is $6$-significant.
One can take $m=[\MG:\bG]$.
\endtheorem

\proof
According to \autoref{conjugate}, any
submodule as in the statement contains vectors conjugate to
$-\be_1+\be_2$ and $t\be_1+\be_2$. Hence, the universal
subgroup has both $2$- and $3$-torsion and, in view of
\autoref{G(N)} and~\eqref{eq.Euler}, it is $6$-significant.

For the converse, consider a $6$-significant subgroup~$G$ and let
$m=[\MG:G]$. Due to Lemmas~\ref{significant} and~\ref{N=6.kernel},
$G\supset\MG''_m$ and $G/\MG''_m$ acts
faithfully on $\AM'_m$. Up to conjugation, one can assume that
$\Gs_2\in G$; then $G/\MG''_m$ has the form
\[*
G=\biggl\{\mat[1,a,0,(-t)^s]\biggm|s\in\CG6,\ a\in I\biggr\}
\]
for some ideal $I\subset\Lambda'_m$, and it is clear that $G$ is the
universal subgroup corresponding to the submodule
$\CV=\Lambda'_m\be_2+I\bv\subset\AM'_m$.
\endproof

\example
According to Theorems~\ref{N=6.BG} and~\ref{th.significant},
there are infinitely many conjugacy classes of $6$-significant
subgroups: they can be classified by the
proper ideals $I\subset\Lambda'$ with the property that
$m\Lambda'\subset I$ for some $m$ prime to~$6$.
Using the list in~\cite{Cummins.Pauli}, one can see that
none of them is a congruence subgroup.

\figure
\def\tl#1#2#3{\relax\hbox to0pt{\hss(#1)\enspace
$#3$\hss}}
\centerline{\vbox{\halign{\hss#\hss&&\qquad\qquad\hss#\hss\cr
 \cpic{N6p7}&\cpic{N6p13}\cr
 \noalign{\medskip}
 \tl a{}{p=7}&\tl b{}{p=13}\cr}}}
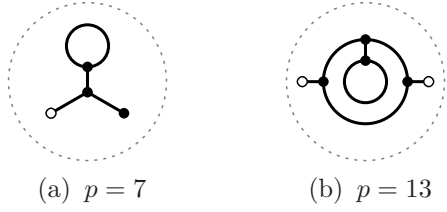
\caption{Examples of~$G_\xi$, $N=6$}\label{fig.N=6}
\endfigure

Two examples
of $6$-significant subgroups
are shown in \autoref{fig.N=6}.
(For the notation, see \autoref{th.N=6} below.)
In each case, the skeleton shown in the figure
and its mirror image correspond
to the two distinct cubic roots of unity in~$\FF{p}$.
\endexample

\subsection{Subgroups with nontrivial type specification}
In conclusion, we show that the only other source of nontrivial
modules $\AM'_G$ is the subgroup $(2A^0)\minus$, see
\autoref{th.MG2}\loccit{MG2.p=0}.

\theorem\label{N=6.notBG}
For a subgroup $G\not\subset\BG3$, denote $\AM'_G=\AM_G/\cp_3$ and
assume that $6\AM'_G\ne0$. Then
\roster
\item\label{N=6.0.BG}
$\IM_G\subset\Lambda(-t\be_1+\be_2)\bmod\cp_3$ and
$G\subset(2A^0)\minus$.
\endroster
If $G$ is $\SS$-transitive and of genus zero,
then $\IM_G=\Lambda(-t\be_1+\be_2)\bmod\cp_3$.
\endtheorem

\proof
Denote by~$\IM'_G$ the image of~$\IM_G$ in $\AM'$. For
any integer $s\ne0\bmod6$,
one has $\Lambda'(t^s-1)\supset6\Lambda'$, and
it follows from~\eqref{eq.N=6} that, whenever the type
specification is not trivial modulo~$6$, there is an inclusion
$6\Lambda'\bv\subset\IM'_G$.
The induced $\BG3$-action on
$\AM'/\bv$ is $\Gb\:h\mapsto(-t)^{\bdeg\Gb}h$; hence, as above,
$\IM'_G$ is \emph{not} a submodule of $6\AM'$ if and only if
$\bG\subset2A^0=\Ker({\bdeg}\bmod2)$ and the type
specification is~$-{\bdeg}\bmod6$, \ie, $G\subset(2A^0)\minus$.

For~$G=(2A^0)\minus$, one does have $\IM'_G=\Lambda'\bv$, and for any
subgroup $G'\subset G$ not contained in~$\BG3$, still
$\IM_{G'}\supset6\Lambda'\bv$. Tensoring the module~$\AM'_{G'}$
with~$\FF2$
or~$\FF3$ and using Theorems~\ref{th.N=3} and~\ref{th.N=2},
one concludes that, if $G'$ is $\SS$-transitive and
of genus zero, the quotient $\Lambda'\bv/\IM_{G'}$ cannot have
$2$- or $3$-torsion.
(Note in addition that the intersection
$(2A^0)\minus\cap\BG3=(\MG')\plus$ is of genus one, see
\autoref{rem.MG}.)
\endproof

To summarize the results obtained in this section,
we restate a few consequences of Theorems~\ref{N=6.BG}
and~\ref{N=6.notBG} in terms of specializations $\IM_G(\xi)$.

\theorem\label{th.N=6}
For an $\SS$-transitive subgroup $G\subset\Bu3$ of genus zero,
assume that the extended Alexander polynomial $\bDelta_{G,p}$ has
a root $\xi\in\bbbk\supset\Bbbk_p$, $p\ne2,3$, with $\ord(-\xi)=6$.
Then one has one of the following two cases\rom:
\roster
\item\label{N=6.0}
$\IM_G=\Lambda(-t\be_1+\be_2)\bmod\cp_3$ and
$G\prec(2A^0)\minus$\rom;
\item\label{N=6.p}
$p\ge5$, $\IM_G(\xi)\sim\bbbk\be_2$, and
$G\prec(G_\xi)\plus\subset\BG3$,
where $G_\xi\subset\MG$ is a certain subgroup of
index $p^{\deg\mp_\xi}$.
\endroster
Cases~\loccit{N=6.0} and~\loccit{N=6.p} are mutually exclusive.
In Case~\loccit{N=6.p}, any \emph{finite} number of \emph{distinct} pairs
$(p,\mp_\xi)$ can appear in the Alexander module of a particular
group.
\endtheorem

\proof
Cases~\loccit{N=6.0} and~\loccit{N=6.p} are given by
Theorems~\ref{N=6.notBG} and~\ref{N=6.BG}, respectively; they are
mutually exclusive due to \autoref{rem.MG}.
In Case~\loccit{N=6.p},
any \emph{finite}
number of distinct primes~$p_i\ge5$ can be `mixed' in the module
$\AM_G=\AM'_m/\Lambda'_m\be_2$, where $m=\prod_ip_i^2$,
see \autoref{N=6.BG}.
\endproof

\let\.\DOTaccent
\def\cprime{$'$}
\bibliographystyle{amsplain}
\bibliography{degt}

\providecommand{\bysame}{\leavevmode\hbox to3em{\hrulefill}\thinspace}
\providecommand{\MR}{\relax\ifhmode\unskip\space\fi MR }
\providecommand{\MRhref}[2]{%
  \href{http://www.ams.org/mathscinet-getitem?mr=#1}{#2}
}
\providecommand{\href}[2]{#2}
\begin{thebibliography}{10}

\bibitem{Artin}
E.~Artin, \emph{Theory of braids}, Ann. of Math. (2) \textbf{48} (1947),
  101--126. \MR{0019087 (8,367a)}

\bibitem{Birch}
Bryan Birch, \emph{Noncongruence subgroups, covers and drawings}, The
  {G}rothendieck theory of dessins d'enfants ({L}uminy, 1993), London Math.
  Soc. Lecture Note Ser., vol. 200, Cambridge Univ. Press, Cambridge, 1994,
  pp.~25--46. \MR{1305392 (95k:11055)}

\bibitem{Bogomolov}
Fedor Bogomolov and Yuri Tschinkel, \emph{Monodromy of elliptic surfaces},
  Galois groups and fundamental groups, Math. Sci. Res. Inst. Publ., vol.~41,
  Cambridge Univ. Press, Cambridge, 2003, pp.~167--181. \MR{2012216
  (2004j:14045)}

\bibitem{Burau}
Werner Burau, \emph{{\"U}ber {Z}opfgruppen und gleichsinnig verdrillte
  {V}erkettungen}, Abh. Math. Sem. Hamburg \textbf{11} (1936), 179--186.

\bibitem{Cummins.Pauli}
C.~J. Cummins and S.~Pauli, \emph{Congruence subgroups of {${\rm PSL}(2,{\Bbb
  Z})$} of genus less than or equal to 24}, Experiment. Math. \textbf{12}
  (2003), no.~2, 243--255. \MR{2016709 (2004i:11037)}

\bibitem{degt:trigonal}
Alex Degtyarev, \emph{Dihedral coverings of trigonal curves}, to appear,
  \verb+arXiv:1005.1038+.

\bibitem{degt:poly}
\bysame, \emph{Alexander polynomial of a curve of degree six}, J. Knot Theory
  Ramifications \textbf{3} (1994), no.~4, 439--454. \MR{1304394 (95h:32042)}

\bibitem{degt:divide}
\bysame, \emph{A divisibility theorem for the {A}lexander polynomial of a plane
  algebraic curve}, Zap. Nauchn. Sem. S.-Peterburg. Otdel. Mat. Inst. Steklov.
  (POMI) \textbf{280} (2001), no.~Geom. i Topol. 7, 146--156, 300, English
  translation: J. Math. Sci. (N. Y.) 119 (2004), no. 2, 205--210. \MR{1879260
  (2002j:14037)}

\bibitem{degt:Oka}
\bysame, \emph{Oka's conjecture on irreducible plane sextics}, J. Lond. Math.
  Soc. (2) \textbf{78} (2008), no.~2, 329--351. \MR{2439628 (2009f:14054)}

\bibitem{degt:dessin}
\bysame, \emph{The fundamental group of a generalized trigonal curve}, Osaka J.
  Math. \textbf{48} (2011), no.~3, 749--782.

\bibitem{degt:monodromy}
\bysame, \emph{{H}urwitz equivalence of braid monodromies and extremal elliptic
  surfaces}, Proc. Lond. Math. Soc. (3) \textbf{103} (2011), 1083--1120.

\bibitem{Esnault}
H{\'e}l{\`e}ne Esnault, \emph{Fibre de {M}ilnor d'un c\^one sur une courbe
  plane singuli\`ere}, Invent. Math. \textbf{68} (1982), no.~3, 477--496.
  \MR{669426 (84a:14003)}

\bibitem{Kulkarni}
Ravi~S. Kulkarni, \emph{An arithmetic-geometric method in the study of the
  subgroups of the modular group}, Amer. J. Math. \textbf{113} (1991), no.~6,
  1053--1133. \MR{1137534 (92i:11046)}

\bibitem{Libgober:Duke}
Anatoly Libgober, \emph{Alexander polynomial of plane algebraic curves and
  cyclic multiple planes}, Duke Math. J. \textbf{49} (1982), no.~4, 833--851.
  \MR{683005 (84g:14030)}

\bibitem{Libgober:invariants}
\bysame, \emph{Alexander invariants of plane algebraic curves}, Singularities,
  {P}art 2 ({A}rcata, {C}alif., 1981), Proc. Sympos. Pure Math., vol.~40, Amer.
  Math. Soc., Providence, RI, 1983, pp.~135--143. \MR{713242 (85h:14017)}

\bibitem{Libgober:modules}
\bysame, \emph{Alexander modules of plane algebraic curves}, Low-dimensional
  topology ({S}an {F}rancisco, {C}alif., 1981), Contemp. Math., vol.~20, Amer.
  Math. Soc., Providence, R.I., 1983, pp.~231--247. \MR{718145 (84k:57002)}

\bibitem{Libgober:Burau}
\bysame, \emph{Invariants of plane algebraic curves via representations of the
  braid groups}, Invent. Math. \textbf{95} (1989), no.~1, 25--30. \MR{969412
  (90a:14038)}

\bibitem{Libgober:survey2}
\bysame, \emph{Problems in topology of the complements to plane singular
  curves}, Singularities in geometry and topology, World Sci. Publ.,
  Hackensack, NJ, 2007, pp.~370--387. \MR{2311493 (2008c:14044)}

\bibitem{Loeser.Vaquie}
F.~Loeser and M.~Vaqui{\'e}, \emph{Le polyn\^ome d'{A}lexander d'une courbe
  plane projective}, Topology \textbf{29} (1990), no.~2, 163--173. \MR{1056267
  (91d:32053)}

\bibitem{Oka:survey}
Mutsuo Oka, \emph{A survey on {A}lexander polynomials of plane curves},
  Singularit\'es {F}ranco-{J}aponaises, S\'emin. Congr., vol.~10, Soc. Math.
  France, Paris, 2005, pp.~209--232. \MR{2145956 (2006m:14038)}

\bibitem{vanKampen}
E.~R. van Kampen, \emph{On the fundamental group of an algebraic curve}, Amer.
  J. Math. \textbf{55} (1933), 255--260.

\bibitem{Zariski:group}
Oscar Zariski, \emph{On the {P}roblem of {E}xistence of {A}lgebraic {F}unctions
  of {T}wo {V}ariables {P}ossessing a {G}iven {B}ranch {C}urve}, Amer. J. Math.
  \textbf{51} (1929), no.~2, 305--328. \MR{1506719}

\end{thebibliography}

\end{document}